\documentclass[preprint,3p,times]{elsarticle}

\usepackage{graphicx}
\usepackage{epsfig}

\usepackage{amssymb}
\usepackage{amsthm}
\usepackage{amsbsy}
\usepackage{amsmath}
\usepackage{amsfonts}
\usepackage{dsfont} 
\usepackage{rotating}
\usepackage{color}

\usepackage{cancel}	 
\usepackage{pgf,pgfplots}

\usepackage{tikz}
\usetikzlibrary{external}
\pgfplotsset{compat=1.15}
\usetikzlibrary{decorations.fractals}
\usetikzlibrary{decorations.pathmorphing,backgrounds}
\usetikzlibrary{decorations.shapes}
\usetikzlibrary{decorations.footprints}
\usetikzlibrary{shapes,arrows,positioning,calc}
\tikzset{paint/.style={ draw=#1!50!black, fill=#1!50 },
    decorate with/.style=
    {decorate,decoration={shape backgrounds,shape=#1,shape size=2mm}}}
\usepackage{mathrsfs}
\usetikzlibrary{arrows}
\definecolor{pergamena}{rgb}{2.44,2.15,0.67}
\definecolor{azz}{rgb}{0.85,0.90,1.00}
\definecolor{mandarino}{RGB}{237,189,101}
\definecolor{f1}{RGB}{255,204,153}
\definecolor{pera}{RGB}{237,210,61}
\definecolor{nonh}{RGB}{237,229,156}
\definecolor{menta}{RGB}{171,242,193}
\definecolor{darkspringgreen}{rgb}{0.09, 0.45, 0.27}
\definecolor{UniBlue}{RGB}{83,121,170}
\definecolor{darkgreen}{rgb}{0.0, 0.2, 0.13}

\newcommand{\reds}[1]{{#1}}

\DeclareMathOperator{\sech}{sech}

\setcitestyle{square}

\date{November, 2017}

\journal{}

\begin{document}

\begin{frontmatter}



\title{Space-time adaptive ADER-DG schemes for dissipative flows: compressible Navier-Stokes and resistive MHD equations}

\author[UNITN]{Francesco Fambri}
\ead{francesco.fambri@unitn.it}
\author[UNITN]{Michael Dumbser\corref{cor1}}
\ead{michael.dumbser@unitn.it}
\author[UNITN]{Olindo Zanotti}
\ead{olindo.zanotti@unitn.it}

\address[UNITN]{Laboratory of Applied Mathematics, Department of Civil, Environmental and Mechanical Engineering, University of Trento, Via Mesiano 77, I-38123 Trento, Italy}

\cortext[cor1]{Corresponding author}

\begin{abstract}
This paper presents an \emph{arbitrary high-order accurate} ADER Discontinuous Galerkin (DG) method on space-time adaptive meshes (AMR) for the solution of two important families of non-linear 
time dependent partial differential equations for compressible \emph{dissipative flows}: the compressible Navier-Stokes equations and the equations of viscous and resistive magnetohydrodynamics 
 in two and three space-dimensions. 

The work continues a recent series of papers concerning the development and application of a proper \emph{a posteriori} subcell finite volume limiting procedure suitable for discontinuous 
Galerkin methods \cite{Dumbser2014,Zanotti2015c,Zanotti2015d}. It is a well known fact that a major weakness of high order DG methods lies in the difficulty of limiting discontinuous 
solutions, which generate spurious oscillations, namely the so-called 'Gibbs phenomenon'. 
In the present work, a nonlinear stabilization of the scheme is sequentially and locally introduced only for troubled cells on the basis of a novel \textit{a posteriori} 
detection criterion, i.e. the MOOD approach. 
The main benefits of the MOOD paradigm, i.e. the computational robustness even in the presence of strong shocks, are preserved and the numerical diffusion is considerably reduced also for the \emph{limited} cells by resorting to a proper sub-grid. 
In practice the method first produces a so-called \textit{candidate solution} by using a high order accurate \emph{unlimited} DG scheme. Then, a set of numerical and physical detection critera is applied to the candidate solution, namely: positivity of pressure and density, absence of floating point errors and satisfaction of a discrete maximum principle in the sense of polynomials. Furthermore, in those cells where at least one of these critera is violated the computed candidate solution is detected as \emph{troubled} and is \emph{locally} rejected. 
Subsequently, a more reliable numerical solution is \emph{recomputed a posteriori} by employing a more robust but still very accurate ADER-WENO finite volume scheme on the subgrid averages within that troubled cell. Finally, a high order DG polynomial is  reconstructed back from the evolved subcell averages. 

We apply the whole approach for the first time to the equations of compressible gas dynamics and magnetohydrodynamics in the presence of viscosity, thermal conductivity and magnetic resistivity, 
therefore extending our family of adaptive ADER-DG schemes to cases for which the numerical fluxes also depend on the gradient of the state vector.  

The distinguished high-resolution properties of the presented numerical scheme stands out against a wide number of non-trivial test cases both for the compressible Navier-Stokes and the 
viscous and resistive magnetohydrodynamics equations. The present results show clearly that the \emph{shock-capturing capability} of the news schemes are significantly enhanced within a 
cell-by-cell \emph{Adaptive Mesh Refinement} (AMR) implementation together with time accurate local time stepping (LTS).

\end{abstract}


\end{frontmatter}


\newcommand{\be}{\begin{equation}}
\newcommand{\ee}{\end{equation}}
\newcommand{\bdm}{\begin{displaymath}}
\newcommand{\edm}{\end{displaymath}}
\newcommand{\bea}{\begin{eqnarray}}
\newcommand{\eea}{\end{eqnarray}}
\newcommand{\PNM}{P_NP_M}
\newcommand{\halb}{\frac{1}{2}}
\newcommand{\FQi}{\tens{\mathbf{F}}\left(\Qi\right)}
\newcommand{\FQj}{\tens{\mathbf{F}}\left(\Qj\right)}
\newcommand{\FQjj}{\tens{\mathbf{F}}\left(\Qjj\right)}
\newcommand{\nj}{\vec n_j}
\newcommand{\FORCE}{\textnormal{FORCE}}
\newcommand{\GFORCE}{\textnormal{GFORCEN}}
\newcommand{\LF}{\textnormal{LF}'}
\newcommand{\LW}{\textnormal{LW}'}
\newcommand{\WL}{\mathcal{W}_h^-}
\newcommand{\WR}{\mathcal{W}_h^+}
\newcommand{\nur}{\boldsymbol{\nu}^\textbf{r} }
\newcommand{\nuf}{\boldsymbol{\nu}^{\boldsymbol{\phi}} }
\newcommand{\nut}{\boldsymbol{\nu}^{\boldsymbol{\theta}} }
\newcommand{\ar}{\phi_1\rho_1}
\newcommand{\arr}{\phi_2\rho_2}
\newcommand{\ur}{u_1^r}
\newcommand{\uf}{u_1^{\phi}}
\newcommand{\ut}{u_1^{\theta}}
\newcommand{\urr}{u_2^r}
\newcommand{\uff}{u_2^{\phi}}
\newcommand{\utt}{u_2^{\theta}}
\newcommand{\ub}{\textbf{u}_\textbf{1}}
\newcommand{\ubb}{\textbf{u}_\textbf{2}}
\newcommand{\RoeMat}{{\tilde A}_{\Path}^G} 
\renewcommand{\u}{\mathbf{u}}
\renewcommand{\v}{\mathbf{v}}
\newcommand{\q}{\mathbf{q}}
\newcommand{\w}{\mathbf{w}}
\newcommand{\U}{\mathcal{U}}
\newcommand{\Q}{\mathcal{Q}}
\newcommand{\dof}[1]{\hat{\mathbf{#1}}}
\newcommand{\real}{{\rm I\!R}}
\newcommand{\nat}{{\rm I\!N}}
\newcommand{\ndof}{N_{\text{\emph{dof}}}}
\newcommand{\proo}[1]{{\color{red} #1 }} 
\newcommand{\CFL}{\text{CFL}} 
\newcommand{\emm}{}  

\definecolor{qqqqff}{rgb}{0.,0.,1.}
\definecolor{ffqqqq}{rgb}{1.,0.,0.}



\section{Introduction}
\label{introduction}

The partial differential equations considered in this paper, namely the compressible Navier-Stokes (CNS) and the viscous and resistive magnetohydrodynamics (VRMHD) equations, 
can be written in a general form that resembles the standard form of a hyperbolic conservation law, except for the fact that diffusivity enters the PDE by means 
of an extra dependence of the flux tensor on the gradient of the solution, i.e. 
\begin{align}
& \frac{\partial \mathbf{u}}{\partial t} + \nabla \cdot \mathbf{F}(\mathbf{u},\nabla \mathbf{u}) = 0, \label{eq:hyp}\\ 
& \mathbf{u} = \mathbf{u}(\mathbf{x},t), \;\;\; \mathbf{x}\in\Omega\subset{\real}^d, \quad t \in {\real}_0^+, \nonumber 
\end{align}
with $\mathbf{u}=\mathbf{u}(\mathbf{x},t)$ being the vector of conserved variables, $\mathbf{F} = \mathbf{F}(\mathbf{u},\nabla \mathbf{u}) = ( \mathbf{f}, \mathbf{g}, \mathbf{h})$ being the nonlinear flux tensor depending in general on the state $\mathbf{u}$ and on its \emph{gradient}  $\nabla \mathbf{u}$. 

Since the compressible Navier-Stokes equations are a special case of the VRMHD equations in absence of electro-magnetic fields 
($\mathbf{B}=0$), we only report the VRMHD equations in the following. The governing equations, which can be cast into the form given by \eqref{eq:hyp} read  
(see \cite{WarburtonVRMHD,DumbserBalsara}): 
\begin{eqnarray}
	  \frac{\partial}{\partial t} \left( \begin{array}{c} \rho \\ \rho \mathbf{v} \\ \rho E \\ \mathbf{B} \\ \psi \end{array} \right) 
		+ \nabla \cdot \left( \begin{array}{c} \rho \mathbf{v} \\ \rho \mathbf{v} \otimes \mathbf{v} + p \mathbf{I} - \boldsymbol{\sigma} - \boldsymbol{\beta} \\ 
		 \mathbf{v} \cdot \left( (\rho E + p) \mathbf{I} - \boldsymbol{\sigma} - \boldsymbol{\beta} \right) - \kappa \nabla T 
		- \frac{\eta}{4 \pi} \mathbf{B} \cdot \left( \nabla \mathbf{B} - \nabla \mathbf{B}^T \right)  \\ 
		\mathbf{B} \otimes \mathbf{v} - \mathbf{v} \otimes \mathbf{B} - \eta \left( \nabla \mathbf{B} - \nabla \mathbf{B}^T \right) + \psi \mathbf{I}  \\ 
		c_h^2 \mathbf{B} \end{array} \right) = 0, 
		\label{eqn.vrmhd} 
\end{eqnarray} 
with the viscous shear stress tensor of the fluid  
\begin{equation}
   \boldsymbol{\sigma} = \mu \left( \nabla \mathbf{v} + \nabla \mathbf{v}^T - \frac{2}{3} \nabla \cdot \mathbf{v} \right),  
\end{equation} 
and the Maxwell stress tensor that contains the stress due to the electro-magnetic forces 
\begin{equation}
	   \boldsymbol{\beta} = \frac{1}{4 \pi} \left( - \halb \mathbf{B}^2 \, \mathbf{I} + \mathbf{B} \otimes \mathbf{B} \right). 
\end{equation} 
In the above equations, $\rho$ is the fluid density, $\mathbf{v}$ is the velocity vector, $\mathbf{B}$ is the magnetic field and 
$p$ is the fluid pressure. The total energy density $\rho E$, which contains the internal energy density $\rho e$, the kinetic 
energy density $\halb \rho \mathbf{v}^2$ and the magnetic pressure $\frac{1}{8 \pi} \mathbf{B}^2$ is related to the 
fluid pressure $p$ by the ideal gas equation of state (EOS):   
\begin{equation}
   \rho E = \frac{p}{\gamma -1} + \frac{1}{2} \rho \mathbf{v}^2 + \frac{1}{8 \pi} \mathbf{B}^2,  
\end{equation} 
where $\gamma$ is the ratio of specific heats. In the VRMHD system given above, $\mu$ is the kinematic viscosity of the fluid, $\kappa$ is the heat conduction coefficient and $\eta$ is
the electric resistivity of the medium. 
The artificial scalar $\psi$ has been introduced in order to deal with the divergence constraint $\nabla \cdot \mathbf{B}=0$ on the magnetic field (which is always true at the continuous level, but
not necessarily at the discrete level inside a numerical scheme),  
according to the generalized Lagrangian multiplier approach (GLM) of Dedner et al., see \cite{Dedneretal}. Here, $c_h$ is an \textit{artificial} propagation speed. The idea of the GLM 
approach is to transport the divergence errors produced by the numerical scheme outside the computational domain. For numerical methods that enforce the divergence condition rigorously also at the discrete level
in the context of the MHD equations, see the work of Balsara et al. \cite{BalsaraSpicer1999,BalsaraAMR,Balsara2004,balsarahllc2d,ADERdivB}. 

We recall that the Navier-Stokes equations are of general theoretical and practical interest for the description of fluid flow with a wide spectrum of applications, ranging from the field of hydraulics, 
oceanic and atmospheric flow modeling, mantle convection in geophysics, aerospace, mechanical and naval engineering up to the simulation of physiological fluid flows in the human cardiovascular or  respiratory system. 
On the other hand, there are many interesting flows of magnetized fluids (plasmas) which are typically described by the MHD equations, but in which resistivity effects of electromagnetic fields 
are also important, such as in solar flares, in the magnetosphere of neutron stars,  in inertial or magnetic confinement fusion for civil energy production, in plasma actuators for active control 
of boundary layers, but also in plasma thrusters for the propulsion of satellites and small spacecraft, just to mention a few examples.  
The dynamics of most of these systems is well approximated by means of the above PDE system (\ref{eq:hyp}). 

Solving the smallest spatial and temporal scales over long time periods and within large domains requires high order of accuracy in both space and time, in order to produce low numerical 
dissipation and dispersion errors. However, discontinuities and very steep gradients can be generated by the above PDE system \eqref{eq:hyp} after finite times, even when starting from perfectly smooth 
initial conditions, due to the non-linearity of the governing equations.
%
Notoriously, finite-volume (FV) methods have been largely used for solving hyperbolic problems and very robust numerical schemes have been developed from them. 
FV schemes are particularly suitable for problems with strong shock waves and are available also on general unstructured meshes but, regrettably, higher order of accuracy in space 
can be achieved only with a loss of simplicity because of the cumbersome \emph{recovery} or \emph{reconstruction} procedures associated to large stencils required by essentially 
non oscillatory (ENO) or weighted ENO (WENO) schemes \cite{abgrall_eno,Dumbser2007}.  
In addition, the corresponding reconstruction stencil introduces a non-trivial spatial-dependence in the computational domain that can deteriorate the parallel scalability of high order 
FV algorithms. Probably because of these complications, discontinuous Galerkin (DG) methods have become increasingly popular over the last decade.

DG methods date back to 1973, when Reed and Hill \cite{reedhill} 
introduced for the first time this new class of finite element methods for solving the neutron transport equation that allows the \emph{'flux to be discontinuous across triangle interfaces'}. 
Reed and Hill noticed a gain in terms of stability with respect to the classical continuous finite element counterpart, but it was only later in the 90'ies that  
the DG formulation was extended to the general case of nonlinear hyperbolic systems in a well-known series of papers by Cockburn and Shu and collaborators 
\cite{Cockburn1989a,Cockburn1989b,Cockburn1990,CockburnShu98}. A review of DG finite element methods is provided in \cite{cockburn_2000_dg,cockburn_2001_rkd}.

In the DG framework, arbitrary high order of accuracy in space can be directly obtained by increasing the polynomial degree of the basis and test functions, i.e. by increasing the number of degrees 
of freedom ($d.o.f$) per element. On the other hand, a stable high order time integration was typically reached through a TVD Runge-Kutta scheme, see  \cite{shu2},  
thus generating the so-called family of RKDG schemes. 
As an alternative to this approach, in this paper the so called \emph{ADER} strategy is used: the ADER approach was introduced in the finite volume context in a series of papers by Toro and Titarev 
\cite{toro1,toro3,toro4,titarevtoro,Toro:2006a}, where \emph{arbitrary} high order of accuracy in space and time is obtained by means of a fully-discrete one-step formulation of the scheme based on 
the numerical solution of a generalized Riemann problem (GRP) at the element interfaces. For a more detailed overview over the GRP, see 
\cite{Ben-Artzi1984,BenArtzi:2006a,Toro2002,titarevtoro,CastroToro,Montecinos2012,GoetzIske}.  
In the original version of the ADER approach a truncated Taylor expansion of the solution in time was used, combined with a sequential analytical differentiation and substitution of the governing equations that replaces time derivatives with space derivatives, i.e. the well-known Cauchy-Kovalewskaya procedure. For general PDE systems this cumbersome procedure leads to a strictly problem 
dependent algorithm, see e.g. \cite{dumbser_jsc,taube_jsc}, with loss of simplicity when discretizing more complex PDE systems or when increasing the desired order of accuracy. 
Furthermore, Taylor series expansions and the standard Cauchy-Kovalewskaya procedure are not able to deal with stiff source terms. For that reason, an alternative version of the ADER approach has been proposed in \cite{DumbserEnauxToro}, in which the Cauchy-Kowalewsky procedure has been replaced by a \emph{local} \emph{space-time} discontinuous Galerkin \emph{predictor}, which is based on a weak formulation of the governing PDE in space and time. The main features of this new version of the ADER approach are threefold: i) it benefits in terms of generality, since the weak formulation in space-time only requires pointwise evaluations of fluxes and source terms, rather than analytical manipulations of the differential form of the PDE; ii) the numerical integration is performed only locally, 
i.e. \emph{within a single space-time element}, minimizing drastically the stencil dependence between the spatial elements and allowing an almost perfect parallelization; iii) the predictor stage is
locally implicit, which allows the treatment of stiff source terms. Over the years, this new formulation of the ADER approach has been successfully applied in several works, including space discretizations based on finite volume schemes \cite{DumbserZanotti,HidalgoDumbser,DumbserNSE,AMR3DCL,AMR3DNC,Zanotti2015,BalsaraMultiDRS} as well as on DG schemes \cite{dumbser_jsc,QiuDumbserShu,LTS,Zanotti2015c,Zanotti2015d}.  

In DG schemes, spurious oscillations may arise when approximating discontinuities and this is a notorious mathematical issue in signal analysis with the name of 'Gibbs phenomenon', noticed for the first time in \cite{Wilbraham}. From the numerical point of view, the experimental observation of the Gibbs oscillations for higher order methods reflects the content of Godunov's theorem \cite{Godunov59},  which states that there is no \emph{linear and monotone} scheme with better than first order of accuracy.  
As a consequence,  over the years several kind of limiters have been developed to avoid Gibbs phenomenon, by resorting, for instance, to the use of \emph{artificial viscosity} \cite{Hartman_02,Persson_06,Feistauer4,Feistauer5,Feistauer6,Feistauer7}, or to 
filtering \cite{Radice2011}, {\emph{a priori}} WENO/HWENO-based reconstruction procedures \cite{QiuShu1,QiuShu3,Qiu_2004,balsara2007,Zhu_2008,Zhu_13,Luo_2007}, or, to slope and moment limiting 
\cite{cbs1,Krivodonova_2007,Biswas_94,Burbeau_2001,Yang_parameterfree_09,Kuzmin2014}. 
In the last few years a novel approach based on a multi dimensional optimal order detection (MOOD) has been introduced with the work of Clain, Diot and Loub\`ere 
\cite{CDL1,CDL2,CDL3,ADER_MOOD_14,ALEMOOD1} in the context of high order finite-volume schemes for multi-dimensional hyperbolic systems on general unstructured meshes. According to the \emph{MOOD paradigm}, the numerical solution is checked only \emph{a posteriori} for some specified mathematical or physical admissibility criteria and then, only for the detected \emph{troubled cells}, the numerical solution is \emph{locally recomputed} according to a \emph{different} and more \emph{robust} numerical scheme, which is typically more dissipative but more stable. In this way the difficulties concerned with the \emph{a priori} prediction of the location of the \emph{future} troubled zones are completely  bypassed.  Because of the encouraging results obtained with this new strategy, the MOOD paradigm has been recently reinterpreted in the context of the DG framework in \cite{Dumbser2014} 
by incorporating a robust higher order ADER-WENO finite volume scheme at the sub-grid level into an arbitrary high order ADER-DG scheme that is active on the main grid. 
The resulting ADER-DG method supplemented by the \emph{a posteriori} ADER-WENO finite volume limiter has been later extended to space-time adaptive Cartesian meshes (AMR) 
\cite{Zanotti2015c,Zanotti2015d}, leading to an unprecedented resolution of shock waves and discontinuities. In this work an attempt is made to extend this new class of schemes also to 
compressible \emph{dissipative flows}. 

At this point, we also would like to point out the very recent work of Peshkov \& Romenski and collaborators \cite{PeshRom2014,HPRmodel}, where the dissipative effects inside a fluid can be 
successfully described within the more general and unified framework of \textit{first order} symmetric \textit{hyperbolic} thermodynamically compatible systems of Godunov \& Romenski, see  
\cite{God1961,God1972MHD,Rom1998,GodRom2003}. The new unified approach is able to describe at the same time viscous fluids as well as elastic and elasto-plastic 
solid media within a single PDE system. 

The rest of the paper is organized as follows: in Section \ref{srh} we first present the numerical method while in the following Section \ref{Results} the computational results
for a large number of non-trivial test cases are shown. Finally, the paper is rounded-off by some concluding remarks in Section \ref{sec.concl}. 

Throughout this paper we will use the Einstein summation convention implying summation over two repeated indices. 


\section{The numerical scheme}
\label{srh}
%
In this section the ADER-DG scheme with \emph{a posteriori} subcell limiter (SCL) on AMR grids is presented in its fundamental facets. For a more detailed description of the scheme, see \cite{Zanotti2015c,Zanotti2015d}.
	
The spatial domain $\Omega$ is discretized with a total number of $N_E$ Cartesian and non-overlapping elements \footnote{In (\ref{eq:mesh}) $\circ$ denotes the \emph{interior} operator, i.e. only the  boundary surfaces of the elements overlap, not the volumes.}  
$\Omega_i$ 
\begin{align}
\Omega =  \bigcup \limits_{i=1,\ldots N_E} \Omega_i, \;\;\; \bigcup \limits_{i\neq j; \;\;i,j=1,\ldots N_E} \Omega^{\circ}_i \cap   \Omega^{\circ}_j =  \varnothing  \label{eq:mesh}
\end{align}
over which we provide the weak formulation of the governing equations (\ref{eq:hyp}), namely 
\begin{align}
\int\limits_{\Omega_i\times T_{n+1}} \phi_k \left( \frac{\partial \mathbf{u}}{\partial t} + \nabla \cdot \mathbf{F}(\mathbf{u},\nabla \mathbf{u})\right)  \,d\mathbf{x}dt = 0. \hspace{2cm}
   i=1,2,\ldots N_E, \;\;\; n \in \nat^+_0 \label{eq:weakPDE}\,.
\end{align}
Here $T_{n+1}=[t^n,t^{n+1}]$ is the current time interval, while $\phi_k\in \U_h^N$ is a generic piece-wise polynomial test-function belonging to the vectorial space  $\U_h^N$ of piecewise polynomials  defined over $\Omega$ and of maximum degree $N\geq 0$, whose discontinuities lie along the element interfaces $\partial \Omega_i$, $i=1$,$2$, $\ldots$ $N_E$. 
As basis and test functions $\phi_k$ we use the set of Lagrange interpolation polynomials of maximum degree $N$ over $\Omega_i$ passing through the Gauss-Legendre quadrature points of the element $\Omega_i$. After integration by parts of the divergence term, equation (\ref{eq:weakPDE}) becomes 
\begin{align}
\int\limits_{\Omega_i\times T_{n+1}} \phi_k  \frac{\partial \mathbf{u}}{\partial t}\,d\mathbf{x}dt + \int\limits_{\partial \Omega_i \times T_{n+1}} \phi_k  \mathbf{F}(\mathbf{u},\nabla \mathbf{u}) \cdot \mathbf{n}   \, dS dt  -  \int\limits_{\Omega_i \times T_{n+1}} \nabla \phi_k  \cdot \mathbf{F}(\mathbf{u},\nabla \mathbf{u})  \,d\mathbf{x}dt= 0, \label{eq:DGPDE}
\end{align}
Notice the total dimension of the chosen space of solutions is $\text{\emph{dim}}(\U_h^N) = N_E \cdot \ndof  $, having $ \ndof =(N +1 )^d$ degrees of freedom ($d.o.f$) for each spatial element $\Omega_i$.
After integrating in time the first term and restricting the space of the solutions to the set of piecewise polynomials $\mathbf{u}_h (\mathbf{x}, t) \in \U_h^N$, i.e. 
\begin{align}
\mathbf{u}_h(\mathbf{x},t) = \phi_k (\mathbf{x})\; \hat{\mathbf{u}}_k (t) 
\end{align}
the following higher order accurate ADER-DG scheme is obtained for the expansion coefficients $\hat{\mathbf{u}}_k^n = \hat{\mathbf{u}}_k(t^n)$: 
\begin{align}
\left( \int\limits_{\Omega_i} \phi_k \phi_l \, d\mathbf{x}\right) \left( \hat{\mathbf{u}}_l^{n+1} - \hat{\mathbf{u}}_l^{n} \right)
+ \int\limits_{\partial \Omega_i \times T_{n+1}} \!\!\!\!\!\! \phi_k  \mathcal{G}\left(\q_h^-, \nabla \q_h^- ;  \q_h^+,\nabla  \q_h^+ \right)  \cdot \mathbf{n}  \, dS \, dt  \; 
- \int\limits_{\Omega_i \times T_{n+1}} \!\!\! \nabla \phi_k  \cdot \mathbf{F}(\q_h,\nabla \q_h)  \,d\mathbf{x}dt= 0, 
\label{eq:ADER-DG}
\end{align}
where a so-called local space-time predictor solution $\mathbf{q}_h(\mathbf{x},t)$ has been introduced and the jumps at the element boundaries are resolved by the (approximate) solution of a Riemann problem at  the element interfaces. In \eqref{eq:ADER-DG} above, the Riemann solver (numerical flux function) is denoted by the symbol $\mathcal{G}\left(\q_h^-, \nabla \q_h^- ;  \q_h^+,\nabla  \q_h^+ \right)$,  depending on a left pair of state $\q_h^-$ and gradient $\nabla \q_h^-$ taken from within the element $\Omega_i$, and a right pair of state $\q_h^+$ and gradient $\nabla \q_h^+$ computed from the adjacent neighbor element, respectively. It has to be noted that even for parabolic equations, an appropriate numerical flux function can be obtained by the solution of a generalized Riemann problem,
see the work of Gassner et al. \cite{MunzDiffusionFlux}, which has also been adopted in \cite{DumbserNSE,DumbserBalsara,HidalgoDumbser}. 
 For the numerical simulations presented in this paper,  $\mathcal{G}$ has been chosen to be a classical and very simple Rusanov-type (local Lax-Friedrichs - LLF) Riemann solver \cite{Rusanov1961a}, 
which has been suitably adapted to account for both hyperbolic and \emph{parabolic} terms, see \cite{DumbserNSE,HidalgoDumbser}:
\begin{equation}
    \mathcal{G}\left(\q_h^-, \nabla \q_h^- ;  \q_h^+,\nabla  \q_h^+ \right) \cdot \mathbf{n} = 
		\frac{1}{2} \left( \mathbf{F}(\q_h^+,\nabla \q_h^+) + \mathbf{F}(\q_h^-,\nabla \q_h^-)  \right) - \frac{1}{2} s_{\max} \left( \q_h^+ - \q_h^- \right), 
\end{equation} 
with 
\begin{equation} 
  s_{\max} = \max{ \left( |\lambda_c(\q_h^-)|, |\lambda_c(\q_h^+)| \right) } + 2 \eta \max{ \left( |\lambda_v(\q_h^-)|, |\lambda_v(\q_h^+)| \right) },  
		\qquad \textnormal{ and } \qquad \eta = \frac{N+1}{h},  
\end{equation} 
where $N$ is the polynomial approximation degree and $h$ is a characteristic length scale of the elements. The $\lambda_c$ denote the eigenvalues of the 
convective (hyperbolic) part of the PDE, i.e. the eigenvalues of the matrix $ \left( \partial \mathbf{F} / \partial \mathbf{u} \right) \cdot \mathbf{n}$, 
while the $\lambda_v$ are the eigenvalues of the parabolic part of the PDE, i.e. those of the matrix 
$ \left( \partial \mathbf{F} / \partial (\nabla \mathbf{u} \cdot \mathbf{n}) \right) \cdot \mathbf{n}$.  
Assuming the space-time predictor $\q_h$ is a polynomial known up to order $(N+1)$ in space and time, see the next paragraphs for the details, then the integrals in (\ref{eq:ADER-DG}) can be
computed \textit{exactly} and the scheme \eqref{eq:ADER-DG} yields an \emph{explicit} and \textit{fully-discrete} one-step formula for the computation of the unknowns at the new time level 
$\hat{\mathbf{u}}_l^{n+1}$. 
For smooth-solutions, the scheme \eqref{eq:ADER-DG} is of order $(N+1)$, see \cite{DumbserNSE}, in principle for \emph{any} integer $N\in\nat^+_0$. On the other hand, a severe time step restriction 
is the curse of all known explicit DG discretizations, i.e. a $\CFL$-type time step restriction of the type 
\begin{align}
\Delta t < \CFL \frac{h_{\text{min}} }{d \left(2N+1\right)}  \left[ \lambda_c^{\text{max}} + \lambda_{v}^{\text{max}} \frac{2(2N+1)}{h_{\text{min}} }\right]^{-1}, \label{eq:CFL}
\end{align}
with the minimum mesh size $h_{\min}$ and $\CFL < 1$. Condition (\ref{eq:CFL}) provides a dependence of the maximum admissible numerical time step $\Delta t$ on the degree $N$ of the polynomial 
basis, the number of space-dimensions $d$, the minimum mesh size given by the insphere diameter $h_{\text{min}}$, the maximum \emph{hyperbolic} signal velocity 
$\lambda_c^{\text{max}}$ and the \emph{parabolic} penalty $\lambda_v^{\text{max}}$ (see \cite{stedg1,stedg2,MunzDiffusionFlux,DumbserBalsara}). 

Equation (\ref{eq:ADER-DG}) is the elementary equation for the time-evolution of the presented ADER-DG-$\mathbb{P}_N$ method. In the following paragraphs the aforementioned local 
space-time predictor $\q_h(\mathbf{x},t)$ and the ADER-WENO subcell limiter, coupled within the space-time AMR framework, are briefly discussed. More details are available in the 
work of \cite{AMR3DCL,AMR3DNC,Zanotti2015,Dumbser2014,Zanotti2015c,Zanotti2015d}. Concerning alternative subcell limiter approaches of the DG method, the reader is referred to 
\cite{Sonntag,CasoniHuerta1,CasoniHuerta2,Fechter1,MeisterOrtleb}. 

\subsection{The element local space-time DG predictor}
\label{ssec:ADERDG}
	
A direct computation of the integral of the non-linear fluxes in equation (\ref{eq:DGPDE}) is subordinate to the knowledge of the physical variables $\u_h$ for any time $t\in T_{n+1}$ along the entire computational domain $\Omega$ or, in other words, to a \emph{fully coupled implicit} solution of the non-linear equation (\ref{eq:DGPDE}) in the coefficients $\dof{u}(t)$  that can become computationally  very demanding. 
Notice that equation system (\ref{eq:ADER-DG}) is already formally conservative, hence it is possible to use a \textit{non-conservative} predictor solution $\q_h$, which can be computed 
\textit{locally} inside each element, without considering any coupling to neighbor elements. In this manner, the resulting computational costs are drastically reduced with respect to the 
original fully coupled system (\ref{eq:DGPDE}). A natural solution to this problem has been presented for the first time in the work of \cite{DumbserEnauxToro} in the context of finite 
volume schemes. 

In this paper we use a \textit{nodal} space-time basis of degree $N$, given by the set of Lagrange interpolation polynomials $\theta_k$ of maximum degree 
$N$ over $\Omega_i\times T_{n+1}$, passing  through the \emph{space-time} Gauss-Legendre quadrature points. 
Since a nodal basis is used, we also expand the nonlinear 
flux tensor as well as the gradient of the solution in the same basis, see \cite{DumbserNSE}. Hence, we have 
\begin{align}
 \q_h(\mathbf{x},t) = \theta_k(\mathbf{x},t) {\hat \q}_k,  \quad 
 \nabla \q_h(\mathbf{x},t) = \nabla \theta_k(\mathbf{x},t) {\hat \q}_k := \theta_k(\mathbf{x},t) {\hat \q}'_k, \quad  
 \mathbf{F}_h(\mathbf{x},t) = \theta_k(\mathbf{x},t) \hat{\mathbf{F}}_k, \,\,\, \textnormal{ with } \,\,\, 
  \hat{\mathbf{F}}_k = \mathbf{F}({\hat \q}_k,{\hat \q}'_k).  
\label{eqn.stdof}
\end{align}
Then, equation (\ref{eq:weakPDE}) reduces to the following element-local system of nonlinear equations for the local space-time predictor polynomials $\q_h(\mathbf{x},t)$: 
\begin{align}
\int\limits_{\Omega_i\times T_{n+1}} \theta_k  \frac{\partial \q_h}{\partial t}\,d\mathbf{x} \, dt  +  \int\limits_{\Omega_i \times T_{n+1}}  \theta_k \nabla \cdot \mathbf{F}(\q_h,\nabla \q_h)  \,d\mathbf{x}dt= 0. \label{eq:predictor}
\end{align}
After integrating the first integral by parts in time, and using the \emph{causality principle} (the current solution depends only on the past, i.e. we use some sort of \textit{upwinding in time}) 
then the following element-local system is obtained: 
\begin{align}
  \int\limits_{\Omega_i} \theta_k(\mathbf{x},t^{n+1}) \q_h(\mathbf{x},t^{n+1})       \, d\mathbf{x}  -  
  \int\limits_{\Omega_i} \theta_k(\mathbf{x},t^{n})   \mathbf{u}_h(\mathbf{x},t^{n}) \, d\mathbf{x}  -  
  \int\limits_{\Omega_i\times T_{n+1}}  \!\!\! \frac{\partial \theta_k}{\partial t} \q_h(\mathbf{x},t) \,d\mathbf{x} \, dt +  
  \int\limits_{\Omega_i \times T_{n+1}} \!\!\! \theta_k \nabla \cdot \mathbf{F}(\q_h,\nabla \q_h)  \,d\mathbf{x}dt= 0,   \label{eq:DOFpredictor} 
\end{align}
which can be solved for the unknown space-time degrees of freedom ${\hat \q}_k$ defined in \eqref{eqn.stdof}. 
Equation (\ref{eq:DOFpredictor}) is solved for each element $\Omega_i$ via a simple iterative method for every $i=1$, $2$, $\ldots$, $N_E$ that has been successfully tested 
with and without stiff or non-stiff source terms in the work of \cite{Dumbser2008,DumbserZanotti}. All the multi-dimensional integrals appearing in the relations above can be 
computed \textit{exactly}, since the solution $\mathbf{q}_h(\mathbf{x},t)$ as well as the fluxes and the gradients are approximated by polynomials of degree $N$ in space and time. 

\subsection{The finite volume sub-cell limiter and adaptive mesh refinement (AMR)}
\label{ssec:SCLAMR}
The high order ADER-DG scheme given by (\ref{eq:ADER-DG}) is an \textit{unlimited} scheme and thus oscillatory in the sense of Godunov. It therefore still 
requires a special treatment for discontinuities. 
Once the local space-time predictor $\q_h(\mathbf{x},t)$ has been obtained from the iterative solution of equation (\ref{eq:DOFpredictor}), as mentioned above, 
then the \emph{candidate solution} $\u^*_h(\mathbf{x},t^{n+1})$ can be directly computed according to equation (\ref{eq:ADER-DG}) in one single step. 
Since the candidate solution $\u^*_h$ may still contain spurious oscillations in the vicinity of steep gradients, underresolved flow features, shock waves or other 
flow discontinuities, nothing can be said about the reliability and about the general physical admissibility of the candidate solution. 
Consequently, a set of physical and numerical admissibility criteria needs to be prescribed and tested. 
A reference point for building  shock-capturing finite-volume schemes is represented by the \emph{discrete maximum principle} (DMP) which is tested on the candidate solution 
accordingly to its \emph{relaxed version} in the sense of polynomials, i.e. in the form 
\begin{equation} 
\label{NAD}
\min \limits_{y\in {\cal{V}}_i} (\textbf{u}_h(\textbf{y},t^n))-\delta \leq \textbf{u}_h^*(\textbf{x},t^{n+1})\leq \max \limits_{\textbf{y}\in {\cal{V}}_i}(\textbf{u}_h(y,t^n))+\delta, 
\qquad \forall \textbf{x}  \in \Omega_i\,,
\end{equation}
where ${\cal{V}}_i$ is the set containing the element  $\Omega_i$ and the respective  Voronoi neighbor elements (neighbors which share a common node with
$\Omega_i$); $\delta$ is chosen to be a \emph{solution-dependent} tolerance given by
\reds{
\begin{equation}
	\delta =\max \left( \delta_0, \epsilon \cdot \left( \max \limits_{y\in {\cal{V}}_i}(\textbf{u}_h(\textbf{y},t^n))- \min \limits_{y\in {\cal{V}}_i}(\textbf{u}_h(\textbf{y},t^n))\right)\, \right),
\end{equation}
}
with $\delta_0=10^{-4}$ and $\epsilon=10^{-3}$, similarly to \cite{Dumbser2014,Zanotti2015c,Zanotti2015d}. The tolerance is added since it is very difficult to compute the global extrema of 
$\mathbf{u}_h(\mathbf{x},t^n)$ in $\Omega_i$. Therefore, we compute an approximation of the extrema by making use of the subgrid representation of the solution, as detailed below. 
Moreover, it is of fundamental importance to check $\u^*_h$ also for a set of \emph{physical} admissibility criteria, e.g. the positivity of pressure and density variables in the case of 
compressible fluid flows. We furthermore check the solution for the presence of floating point errors (NaN). 
Once the numerical and physical admissibility criteria have been tested and whenever a local candidate solution $\u^*_h(\mathbf{x}\in\Omega_i,t^{n+1})$ is detected to be '\emph{troubled}', then  $\u^*_h(\mathbf{x}\in\Omega_i,t^{n+1})$ is directly rejected and the limiter-status of $\Omega_i$ is set to $\beta_i = 1$, meaning the limiter is activated. Then, the older ADER-DG solution $\u_h(\mathbf{x}\in {\Omega}_i,t^{n})$ is projected along a suitable \emph{sub-grid of} $N_s$ \emph{spatial} \emph{sub-cells per space-dimension} within $\Omega_i$, resulting in a \emph{piecewise-constant} representation  of the discrete solution \reds{$\w_h(\mathbf{x}\in {\Omega}_i,t^{n})=\mathcal{P}[\u_h(\mathbf{x}\in {\Omega}_i,t^{n})]$}, $\mathcal{P}$ being a suitable projector operator (see \cite{Dumbser2014,Zanotti2015c,Zanotti2015d}). Then, a new discrete solution is obtained for the subgrid averages by using a more robust ADER-WENO finite volume scheme \cite{AMR3DCL}, generating a new set of piecewise-constant cell averages 
$\w_h(\mathbf{x}\in \Omega_i,t^{n+1})$. The new subcell averages are then directly gathered back to a high order DG polynomial $\u_h(\mathbf{x}\in \Omega_i,t^{n+1})= \mathcal{R}[ \w_h(\mathbf{x}\in \Omega_i,t^{n+1})]$, where $ \mathcal{R}$ is a suitable high order accurate \emph{reconstruction} operator satisfying $ \mathcal{R}  \circ   \mathcal{P} = 1 $ (see \cite{Dumbser2014,Zanotti2015c,Zanotti2015d}).  The high order ADER-WENO method has been shown to be an excellent candidate for the subcell finite volume limiting stage because of its well established capabilities in handling discontinuities, 
together with high-order convergence properties under the time-step constraint 
\begin{align}
\Delta t < \CFL \frac{h_{\text{min}} }{d N_s}  \left[ \lambda_c^{\text{max}} + \lambda_{v}^{\text{max}} \frac{2 N_s}{h_{\text{min}} }\right]^{-1}. \label{eq:CFLweno}
\end{align}
The WENO scheme furthermore does not clip local extrema, in contrast to standard second order TVD schemes. Notice that the local number of sub-cells $N_s$ per space-dimension should be chosen 
$N_s \geq N+1$ in order to preserve the information contained in the available degrees of freedom of the high order polynomial data representation used in the DG scheme. In our simulations, $N_s$ has been chosen to be $N_s=2N+1$, 
thus matching the maximum time-step allowed by the ADER-WENO finite volume scheme \eqref{eq:CFLweno} with the one for the ADER-DG method (\ref{eq:CFL}) .

Just a few words are necessary to briefly introduce the space-time adaptive mesh (AMR) in which the complete numerical scheme is mounted. 
Further details are available in the recent papers of \cite{Zanotti2015c,Zanotti2015d}, where essentially the same AMR technique is used. There exist essentially two different ways of implementing 
an AMR method, and both of them are characterized by \reds{pros} and cons: the first technique is based on the nested structure of \emph{independent overlaying sub-grid 'patches'} (see \cite{Berger-Oliger1984,berger85,Berger-Colella1989}); the second one is the so called \emph{'cell by cell' refinement} and this is the adopted AMR-approach because of its formally very simple tree-type data structure (see \cite{Khokhlov98,AMR3DCL}). In the here-presented 'cell-by-cell' AMR \emph{every single element is recursively refined}, from a coarsest refinement level $\ell_0=0$ to a prescribed finest (maximum) refinement level $\ell_{\text{max}}\in\nat^+_0$, accordingly to a \emph{refinement-estimator function} ${\chi}_{\emm}$ that drives  step by step the choice for recoarsening or refinement. ${\chi}_{\emm}$ is chosen to be
%
\reds{\begin{equation} \label{eq:chi_m}
	{\chi}_{\emm} (\Phi) =  \sqrt{ \frac{\sum_{k,l}{\left(\left. \partial ^2 \Phi \middle/ \partial x_k \partial x_l \right. \right)^2}}{ \sum_{k,l}{ \left[ \left. \Big( \left| \left. \partial \Phi \middle/ \partial x_k\right.\right|_{\text{F}} + \left| \left.\partial \Phi \middle/ \partial x_k \right. \right|_{\text{B}} \Big) \middle/ \Delta x_l \right. + \epsilon \Delta( \Phi, \nabla \Phi, \nabla^2 \Phi )  \right]^2}}  }\,,
\end{equation}
 following \cite{Loehner1987}, and it involves up to the second order derivative of an \emph{indicator function} $\Phi$, chosen to be a mathematical quantity of physical interest that varies point by point in the computational domain, e.g. the local density of the fluid, or the pressure, the vorticity, or an arbitrarily  chosen different function. The summation $\sum_{k,l}$ runs over the available space-directions; the partial  derivative are evaluated forward and backward in space, correspondingly to subscripts F and B, respectively, based on the cell average, in order to sample information also in the vicinity of the local space-element $\Omega_i$; finally, $\Delta( \Phi, \nabla \Phi, \nabla^2 \Phi )$ is a natural majorization of a centered second order derivative. In this work, $\Delta$ takes the form of $\Delta_i = |\Phi|_B + 2 |\Phi|_i + |\Phi|_F$, and it is modulated by the filter-parameter $\epsilon$. This definition helps in preventing unnecessary mesh-refinement in presence of  ripples in the numerical solution. In this work, $\epsilon$ is chosen to be equal to $\epsilon=10^{-3}$, taken constant during the simulations. Throughout this paper we have simply used the fluid density as indicator function, i.e. $\Phi=\rho$. A more sophisticated choice could be used, for example based on the entropy, as recently proposed in the context of CWENO schemes on adaptive meshes, see 
\cite{PS:entropy,SCR:CWENOquadtree,CS:epsweno} for more details. }

\begin{figure}
	\centering
		\includegraphics[trim= 1cm 2cm 3cm 3cm, width=0.5\textwidth]{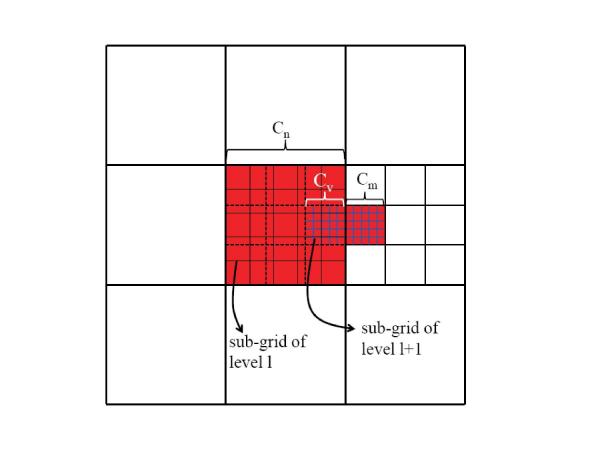}
		\includegraphics[width=0.45\textwidth]{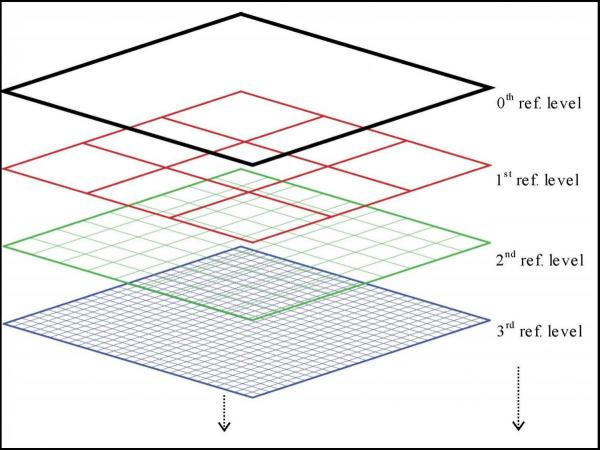}
	\caption{At the left, a simple sketch of the combination of AMR and DG sub-cell reconstruction. The limited elements ($\beta=1$) are highlighted in red, i.e. cell $\mathcal{C}_n$ at the ref. level $\ell$ and cell  $\mathcal{C}_m$ at $\ell+1$. Then, the cell $\mathcal{C}_n$ must project $\textbf{v}_h$ from the original sub-grid of the $\ell$-th ref. level to the sub-grid of level $\ell+1$, within the virtual cell $\mathcal{C}_v$. At the right, the \emph{tree-structure} of the refinement levels  for a single element at the coarsest level $\ell_0$ is shown, \reds{generated after choosing a refinement factor $\mathfrak{r}=3$.} 
(See colored version on-line)}
	\label{fig:AMR}
\end{figure}

\begin{figure} 
\centering
		\includegraphics[width=0.9\textwidth]{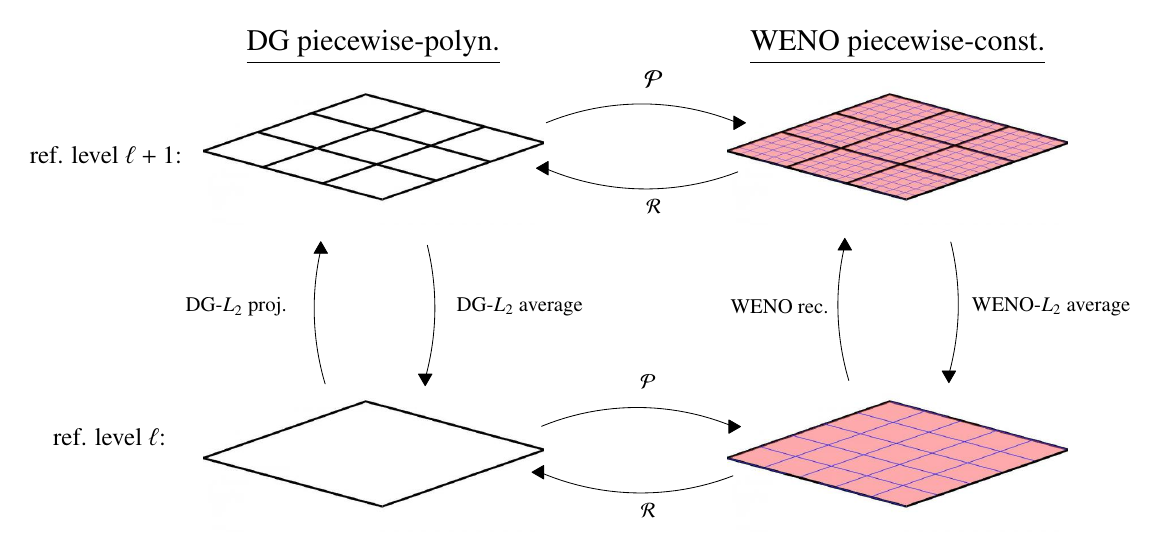} 
\caption{Mapping of the numerical solution between the DG piecewise polynomial  and the WENO piecewise constant spaces, between two different AMR-levels $\ell$ and $\ell+1$.}
\label{fig:AMRmaps}
\end{figure}

A prescribed \emph{refinement factor} $\mathfrak{r}$ indicates the number of sub-element per space-dimension which are generated in a refinement process. \reds{In particular, after a space-element is chosen for refinement, then, 
it is replaced by $\mathfrak{r}^d$ smaller Cartesian space elements (see Fig. \ref{fig:AMR}), and \emph{vice versa} for the recoarsening process. Note that in our code a general refinement  
factor $\mathfrak{r} \geq 2$ can be used. This allows to realize very general space-trees, not only quadtrees or octrees, which would correspond to the choice $\mathfrak{r}=2$.}
Whenever the refinement-process is executed, the refinement-estimator function ${\chi}_{\emm}$ is evaluated in the computational domain and every single element is refined or recoarsened every time ${\chi}_{\emm}$ overpasses the prescribed upper or lower threshold ${\chi}_{\text{ref}}$ and ${\chi}_{\text{rec}}$, respectively. 
It becomes useful, for practical purposes, to assign a \emph{refinement status} $\sigma$ to all the elements at the coarsest level $l_0$ and along the respective finer levels $\ell$ with ${\ell}_0 \leq \ell \leq {\ell}_{\text{max}}$ with the rule
\begin{align}
\sigma_i = \left\{\begin{array}{rcl} 
-1, & & \text{for the so called \emph{virtual parent cells}} \\
0, & & \text{for \emph{active elements}} \\ 
1, & & \text{for the so called \emph{virtual children}} 
\end{array}\right.\;\;\; i=1,\ldots, N_{\text{tot}}.
\end{align}
The \emph{active elements} $\Omega_i$ (i.e. $\sigma=0$) are those \emph{non-overlapping} spatial elements that constitute the current numerical mesh, i.e. satisfying (\ref{eq:hyp}), where the discrete  solution is chosen for being updated following the presented ADER-DG$+$SCL method. 
The \emph{virtual children} belonging to a relative refinement level $\ell(\Omega_i^{\text{Vc}})$ are those spatial elements which \emph{are (spatially) contained within at least  one active element in its } -adjacent coarser- \emph{tree-structure} $\ell = \ell(\Omega_i^{\text{Vc}}) -1$. The numerical solution is updated in time by means of a standard $L_2$ projection for the ADER-DG from the \emph{mother} cell at the $(\ell(\Omega_i^{\text{Vc}})-1)$-th level.
Finally, the \emph{virtual parent cells} $\Omega_i^{\text{Vm}}$ ($\sigma=-1$) belonging to a relative refinement level $\ell(\Omega_i^{\text{Vm}})$ are those spatial elements which \emph{(spatially) contains at least one active element in its adjacent} -finer- \emph{tree-structure} $\ell = \ell(\Omega_i^{\text{Vm}})+1$. In this case, the numerical solution is updated in time by averaging the solution from the \emph{children}-elements (i.e. $\ell=\ell(\Omega_i^{\text{Vm}})+1$).
Within this new computational grid, $N_{\text{tot}}$ is the total number of elements that should be distinguished from the total number of \emph{active} elements $N_{E}$ which appears i n our numerical equations.
These three $\sigma$-status are necessary during the mesh-adaptation stage whenever an active cell is refined or recoarsened, and then, inactivated. Indeed, whenever this is the case, a proper transformation is needed for mapping the numerical solution (limited or unlimited, finer or coarser) from one refinement level to the adjacent one.  
A simple sketch of the transformation-mapping between the discrete solution spaces of the DG polynomials and the WENO subcell averages, between two adjacent refinement levels $\ell$ and $\ell +1$ is shown in Fig. \ref{fig:AMRmaps}.

\reds{The presented AMR framework is known as \emph{h}-refinement. An alternative approach is the so-called \emph{p}-refinement, where the mesh is fixed in time, 
but the polynomial degree of the solution is allowed to vary in every time-step and from one element to another, depending on the desired numerical resolution. 
However, \emph{p}-refinement is only useful in areas of the computational domain where the solution is sufficiently smooth. 
The simplest way for implementing \emph{p}-refinement is introducing a \emph{modal basis} in order to minimize the cost of the projection from the space of piece-wise 
polynomials of degree $N$ and the space of piece-wise polynomials of degree $M \neq N$, see e.g. \cite{Dumbser2007c} for more details on \emph{p}-refinement in the 
context of ADER-DG schemes with local time stepping.  
The main advantage of the nodal basis adopted here is that it gives an immediate knowledge of the value of the solution at the Gaussian quadrature points needed in the weak 
formulation, which makes it computationally very efficient compared to a modal basis, where this information requires additional matrix-vector products. 
However, the extension of this work to \emph{hp}-refinement is in principle feasible, but out of the scope of the present work.}

Moreover, virtual cells allow us to perform polynomial WENO reconstructions along the same refinement level, independently on the effective refinement level of two adjacent active elements (see Fig.  \ref{fig:AMR}). For more details see \cite{Zanotti2015c,Zanotti2015d}, for information about the parallel message passing interface (MPI) implementation of the presented AMR framework see \cite{AMR3DNC,AMR3DCL}. 
It should be noticed that whenever an automatic adaptation of the grid is used, the scheme can in principle handle simultaneously small and large spatial scales. However, due to the CFL condition, also a characteristic time scale is implied by the local mesh spacing. Thus, a proper time-accurate and fully conservative local time stepping (LTS) method has to be adopted in order to use the smaller 
time-steps only for the smaller spatial elements, in favor of performance (see \cite{AMR3DCL}). Finally, it should be mentioned that in our formulation two adjacent active elements are allowed to 
belong to two different refinement levels with the constraint $|\Delta \ell| \leq 1$, i.e. the two elements belong to the same or to an adjacent, finer or coarse, AMR-level. 

\begin{figure} 
		\includegraphics[width=0.9\textwidth]{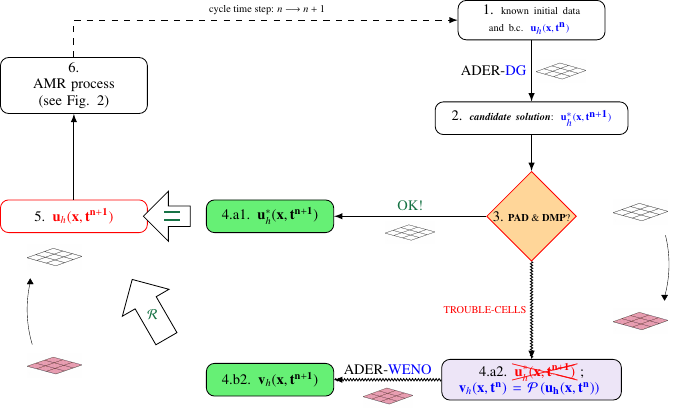}  
\caption{\reds{Simplified flow diagram of the presented numerical method: 1) known the initial solution $u_h$ at time $t^n$ along a given computational AMR mesh; 2) then a \emph{candidate solution} $u_h^*$ is computed through our high-order ADER-DG; 3) then physical admissibility detection, criteria (PAD) and/or mathematical, e.g. the discrete maximum principle (DMP), are checked along the complete computational domain; 4) then two branches are distinguished accordingly to the response of the trouble-cell detection: 4.a1) whenever the \emph{candidate solution} $u_h^*$ is accepted as it is, then the solution $u_h$ is updated according to it; 4.a2) whenever a trouble cell is detected, then the \emph{candidate solution} $u_h^*$ is \emph{directly rejected} only within the trouble cell, and new initial data are generated through simple projection from the DG piecewise-polynomial space to the FV piece-wise constant, i.e. $v_{h}(\bf{x},t^n) = \mathcal{P}\left(u_h(\bf{x},t^n)\right)$; 4.b2) then, accordingly to the FV initial data, a new numerical solution is computed accordingly to a more-robust numerical scheme, e.g. a high-order ADER-WENO finite-volume method; 5) then, a corresponding piece-wise polynomial solution is generated through a reconstruction operator (see \cite{Dumbser2014,Zanotti2015c,Zanotti2015d}); 6) finally, the numerical solution go through our recoarsening/refinement process and the time-step is cycled.}}
\label{fig:MOODflow}
\end{figure}

\reds{A flow diagram depicting the main stages of the final algorithm presented in this section is shown in Fig. \ref{fig:MOODflow}.}


%


\section{Numerical tests}
\label{Results}
\reds{In this section, the accuracy and robustness of our ADER-DG scheme with subcell FV limiting and AMR are tested against a series of non-trivial numerical problems. Up to seven benchmark 
scenarios are chosen with the following criteria: i) an analytical, numerical or experimental reference solution exists;  ii) Mach numbers ranging from medium to high Mach number flow including
strong shock waves, flows with low to high Reynolds regimes are tested; iii) the proper dynamics of physical instabilities and the correct energy dissipation rate can be tested; iv) robustness 
can be tested against shocks and conservation of the physical constraints ($\nabla \cdot \mathbf{B}=0$); v) slip and no-slip wall boundary conditions are considered; vi) the proposed  
\emph{a posteriori} limiting procedure and AMR framework should increase robustness and high resolution properties of the main numerical method without affecting the physics, i.e. no spurious 
mesh effects or unphysical dissipation should be present; vii) a convergence table should be computed against some smooth \emph{analytical} reference solutions.} 

\subsection{Lid-driven cavity flow at low Mach number (M=0.1)}
\label{test:LDC2D}
\begin{figure} 
\centering 
			\includegraphics[width=0.49\textwidth]{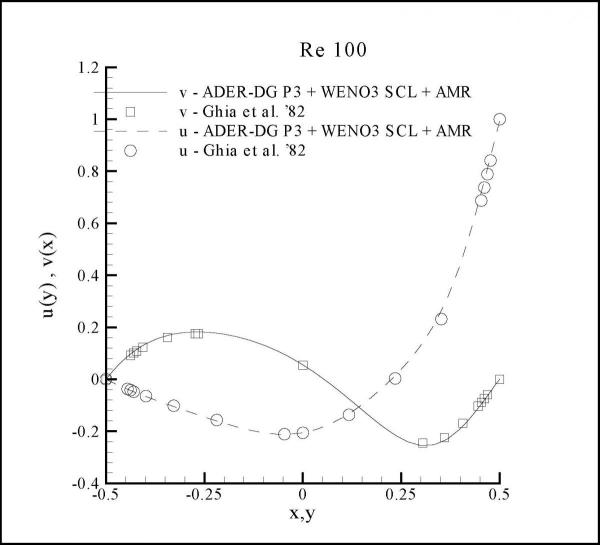} 			\includegraphics[width=0.49\textwidth]{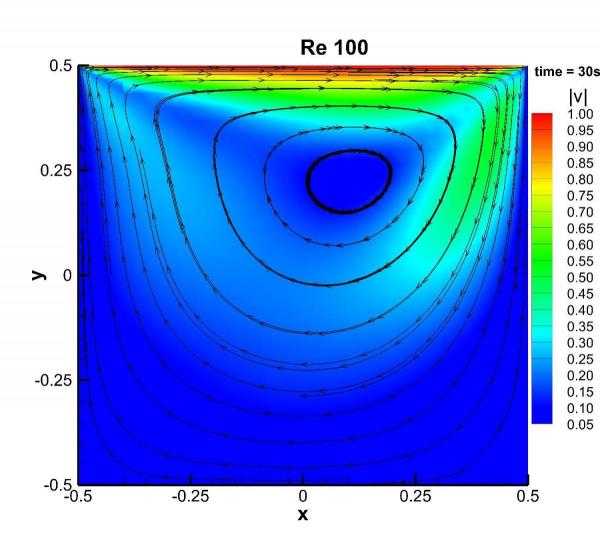}\\
			\includegraphics[width=0.49\textwidth]{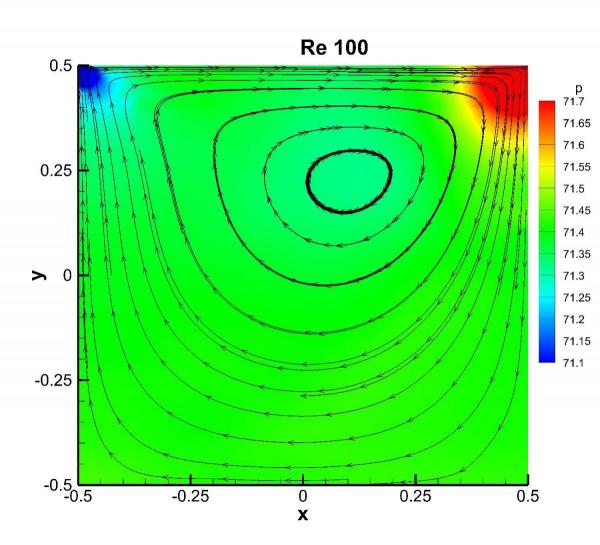}     	\includegraphics[width=0.49\textwidth]{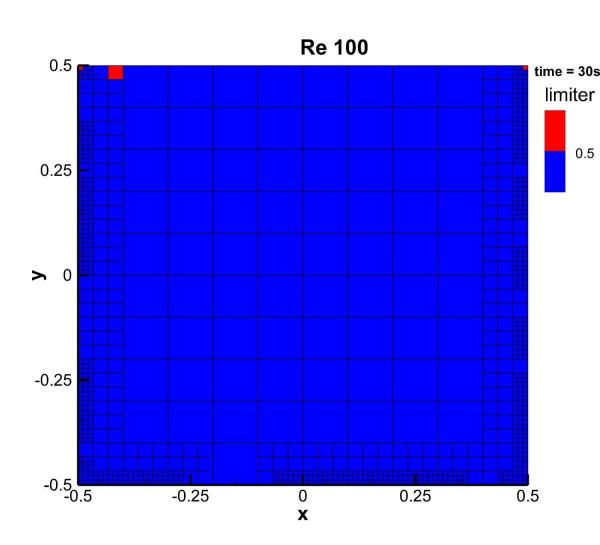} 
\caption{\textbf{Lid-driven cavity.} The numerical solution obtained for the two dimensional lid-driven cavity problem compared with the numerical results of \cite{Ghia1982} at different Reynolds $Re=100$, obtained with our ADER-DG-$\mathbb{P}_3$ method using $10x10$ elements at the coarsest level, up to $\ell_\text{max}=2$ maximum number of refinement levels with  a refine factor $\mathfrak{r}=3$. In the first two rows, from left to right, from the top to the bottom the data-comparison, the magnitude of the velocity field with streamlines, the pressure and the limiter status have been plotted.}\label{fig:LDC2D}
\end{figure}

 The so-called lid-driven cavity flow became a standard benchmark problem for testing numerical methods for the incompressible Navier-Stokes equations, see \cite{Ghia1982}. In this two dimensional test a \emph{nearly} incompressible flow is considered. In a \emph{closed} square cavity the fluid-flow is driven by the moving upper-wall with tangential velocity $u=1$. No-slip boundary conditions are applied at the remaining three walls.  The spatial domain $\Omega=[-1,1]\times[-1,1]$ has been discretized into $10\times 10$ space-elements for the coarsest mesh at level zero; the AMR-framework has been activated accordingly to a refine factor $\mathfrak{r}=3$ and $\ell_\text{max}=2$, the associated maximally refined mesh-level, and the 
magnitude of the velocity as estimator-function for the mesh adaptation. We compare the numerical solution obtained with our ADER-DG-$\mathbb{P}_3$ supplemented with the \emph{a posteriori} WENO$3$ SCL for the compressible Navier-Stokes equations in the low-Mach regime ($M=0.1$) with the reference solution of \cite{Ghia1982} in Fig. \ref{fig:LDC2D}. A very good agreement between computed and reference solution has been obtained, despite the compressibility of the simulated fluid-flow and the non-trivial singularities at the upper corners. 
Notice that the limiter has been needed only next to the flow singularities at the upper corners where, in fact, the solution is a double valued function, i.e. $u=0$ at the side walls and $u=1$ at the moving  upper lid.

\subsection{3D Taylor-Green vortex at low Mach number  (M=0.1)}
\label{test:TGV3D}
\begin{figure} 
\centering 
			\includegraphics[width=0.6\textwidth]{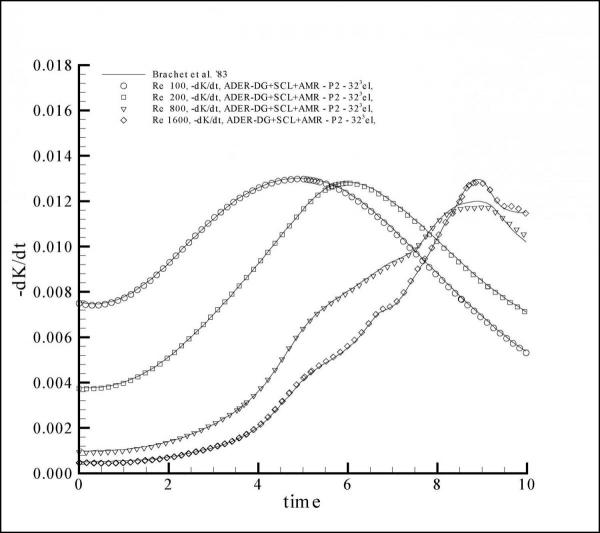} 
			\caption{\textbf{3D Taylor-Green vortex.} Time evolution of the kinetic energy dissipation rate $\epsilon(t)$ obtained with our ADER-DG-$\mathbb{P}_2$ supplemented with the \emph{a posteriori} WENO$3$ SCL	at different Reynolds numbers $Re=100$, $200$ and $800$. The DNS reference solutions of Brachet et al. \cite{Brachet1983} are plotted as continuous lines. Up to $\ell_\text{max}=2$ maximum number of refinement levels with  a refine factor $\mathfrak{r}=2$ are used along the $32^3$ elements of the coarsest grid.} \label{fig:TG3D}
\end{figure}

\begin{figure} 
\centering 
			\includegraphics[width=0.32\textwidth]{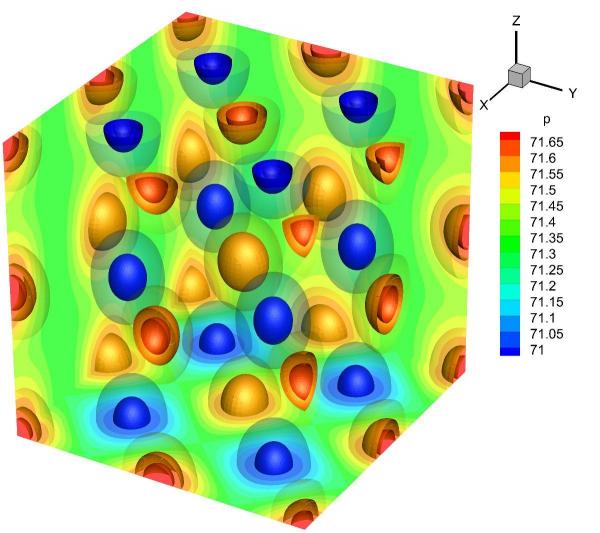} 
			\includegraphics[width=0.32\textwidth]{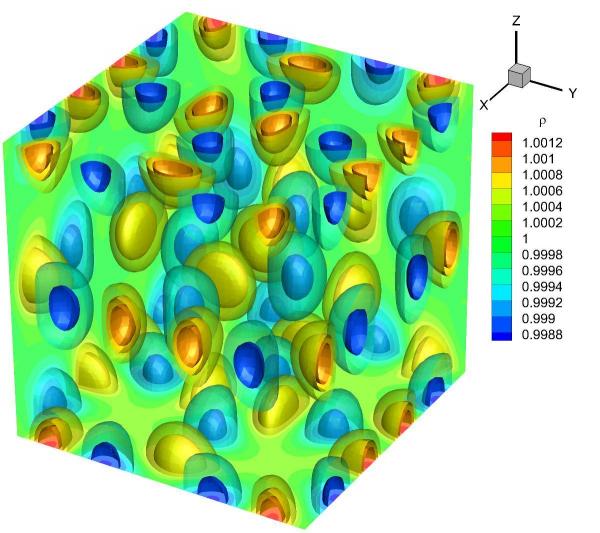} 
			\includegraphics[width=0.32\textwidth]{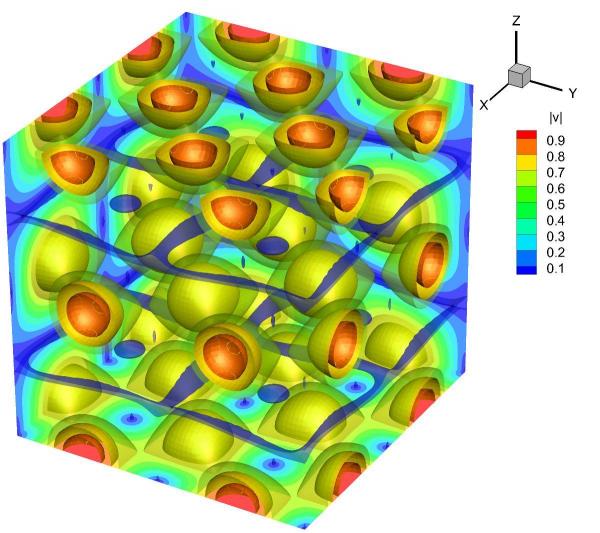} 
			\\ \includegraphics[width=0.32\textwidth]{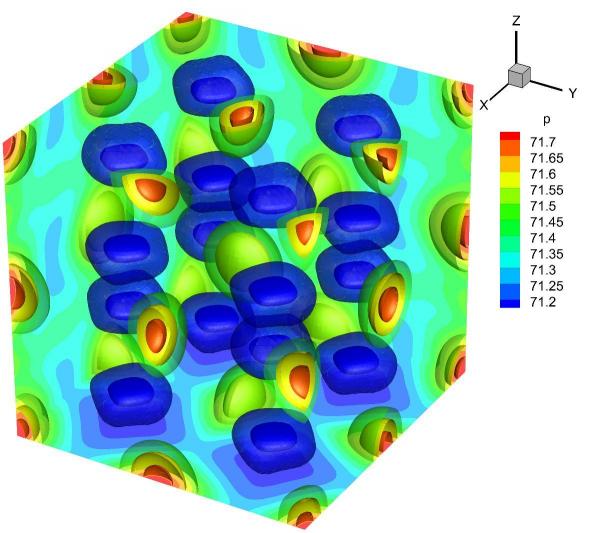}  
			\includegraphics[width=0.32\textwidth]{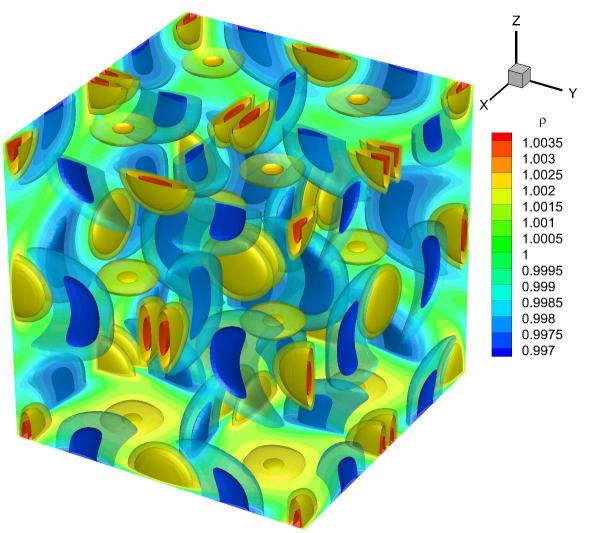} 
			\includegraphics[width=0.32\textwidth]{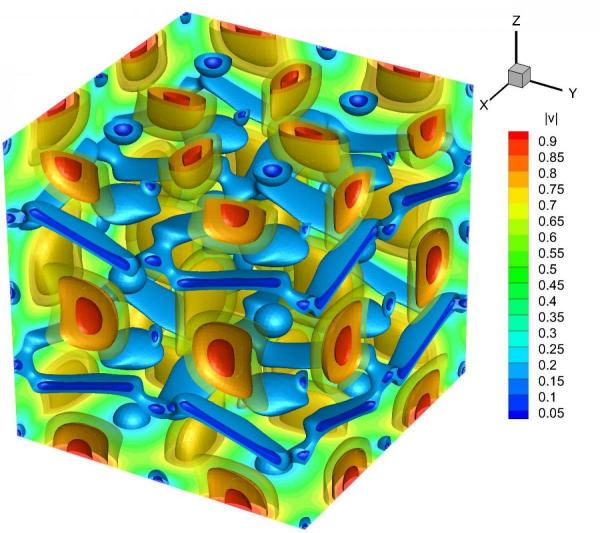}  
			\\ \includegraphics[width=0.32\textwidth]{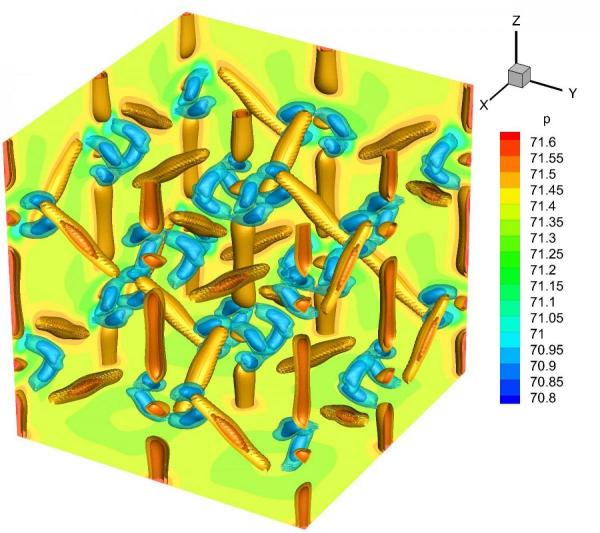} 
			\includegraphics[width=0.32\textwidth]{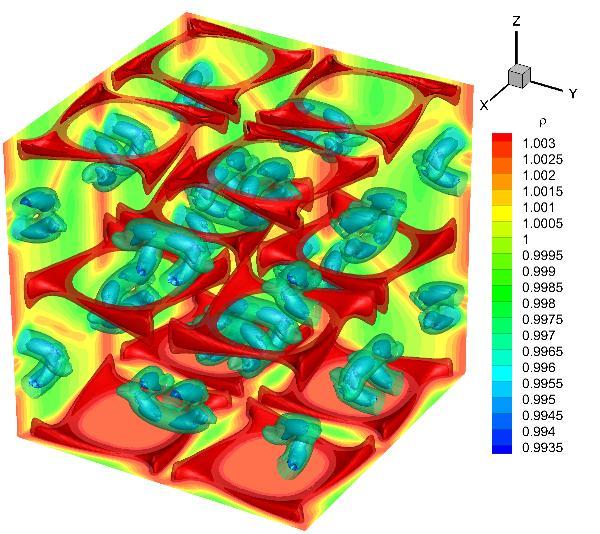} 
			\includegraphics[width=0.32\textwidth]{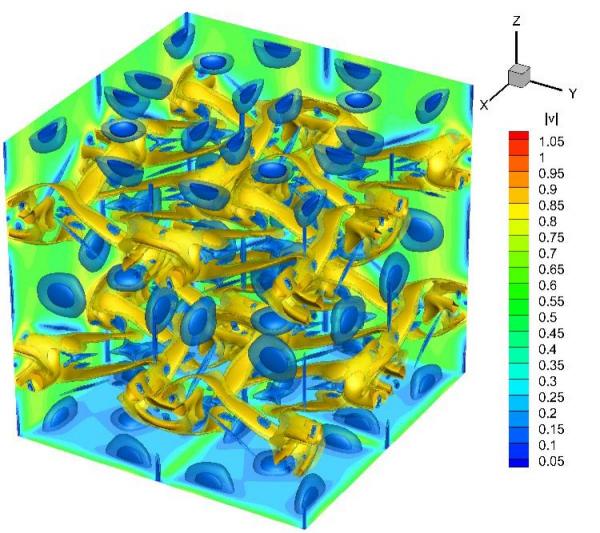}   
			\\ \includegraphics[width=0.32\textwidth]{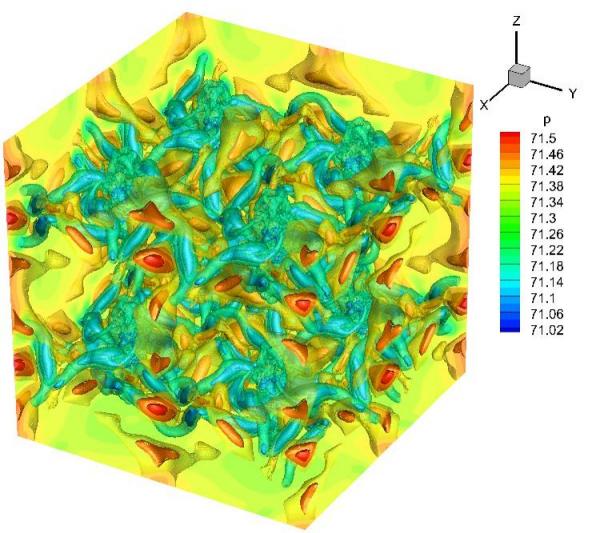}  
			\includegraphics[width=0.32\textwidth]{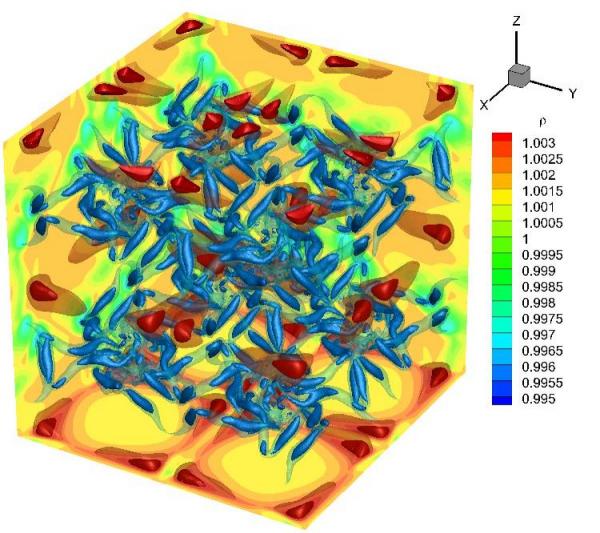} 
			\includegraphics[width=0.32\textwidth]{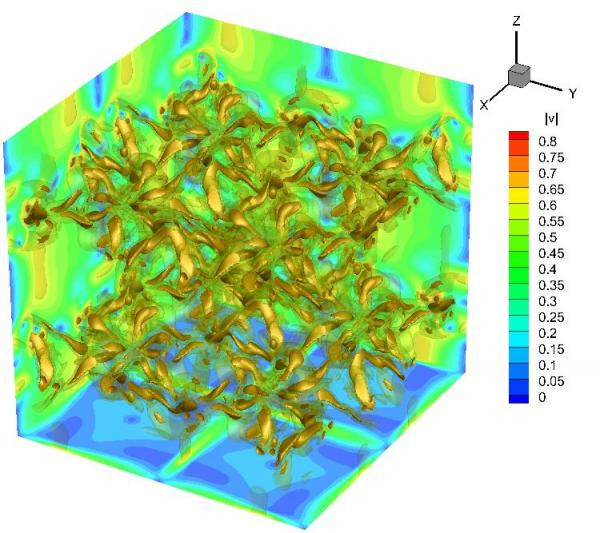}  
			\caption{\textbf{3D Taylor-Green vortex.} Numerical solution for the three dimensional Taylor-Green vortex flow at $Re=800$ computed with our ADER-DG-$\mathbb{P}_{2}$ supplemented by the \emph{a posteriori} SCL using $32^3$ elements on the coarsest level, up to $\ell_\text{max}=2$ maximum number of refinement levels with  a refine factor $\mathfrak{r}=2$.
			The isosurfaces of the pressure (left), the density (center) and the velocity (right) are plotted at 
			times $t=0.5$, $2.0$, $6.0$ and $10.0$ from the top to the bottom, respectively.}\label{fig:TG3D2}
\end{figure}

A very intriguing three-dimensional flow that drives the larger to the smallest physical scales is the turbulence-decaying process that is generated in the Taylor-Green vortex problem. The initial condition of the fluid variables is given by 
\begin{align}
&\rho(x,y,z,0) = 1,\\
&u(x,y,z,0) = \sin(x)\cos(y)\cos(z),\\
&v(x,y,z,0) = -\cos(x)\sin(y)\cos(z),\\
&w(x,y,z,0) = 0,\\
&p(x,y,z,0) = \rho c_0^2/\gamma + \left(\cos(2x)+ \cos(2y) \right)\left(\cos(2z)+2\right)/16.
\end{align}
where $c_0$ is the adiabatic sound speed. The reference solution is widely accepted to be the DNS solution presented by  \cite{Brachet1983} through both a direct spectral 
method (up to $256^3$ modes) and a rigorous power series analysis (up to order $t^{80}$), see also \cite{Morf1980}. 
Periodic boundary conditions are assumed everywhere with respect to the cubic spatial domain $\Omega=[0,2\pi]^3$. 
Fig.  \ref{fig:TG3D} shows the results for the kinetic energy dissipation rate 
\begin{align}
\epsilon(t) = -\frac{\partial K}{\partial t} = - \frac{1}{||\Omega||} \frac{\partial}{\partial t} \int_{\Omega} \frac{1}{2}\rho \mathbf{v}^2 d\mathbf{x}
\end{align}
for different Reynolds numbers $Re\in[100,1600]$ evaluated in the time interval $t\in[0,10]$. A direct comparison with the reference solution of \cite{Brachet1983} shows that an excellent 
agreement has been obtained. Notice that for larger Reynolds numbers, smaller dissipative vortex structures can be generated and, consequently, a higher numerical resolution is needed. 
Since we use a dissipative scheme (due to the Riemann solver), a too low resolution would generate an excess of numerical diffusion. 
The initial condition is the same for all the different test cases, but the time series of the kinetic energy dissipation strongly depends on the chosen Reynolds number. 
At $t=0$ a very smooth solution is initialized, then the diffusive decaying begins slowly. Once the peak of dissipation is reached ($t\sim 4$ for $Re=100$, $t\sim6$ for $Re=200$, $t\sim9$ for $Re=800$ and $Re=1600$) then the kinetic energy dissipation rate decreases asymptotically and inexorably to the trivial stationary solution with $K=0$. For this test the third order $\mathbb{P}_2$ version of 
our ADER-DG scheme supplemented with the third order ADER-WENO$3$ sub-cell limiter has been used. The AMR grid is activated using the $32^3$ elements of the coarsest level zero grid. A refinement 
factor of $\mathfrak{r}=2$ is used and up to $\ell_\text{max}=2$ maximum number of refinement levels are admitted. 
Fig.  \ref{fig:TG3D2} shows the iso-surfaces of pressure, density and velocity at different times $t\in[0,10]$ and gives a better qualitative comprehension of the flow dynamics.

\subsection{Compressible 2D mixing layer} 
\label{test:CML}
\begin{figure} 
\begin{tabular}{l} 
			\includegraphics[width=0.9\textwidth]{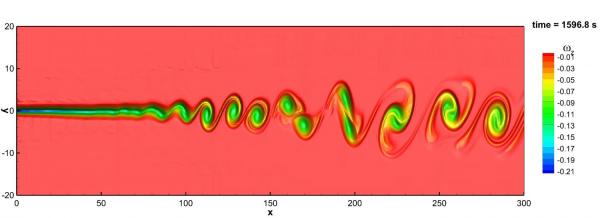}\\
			\includegraphics[width=0.87\textwidth]{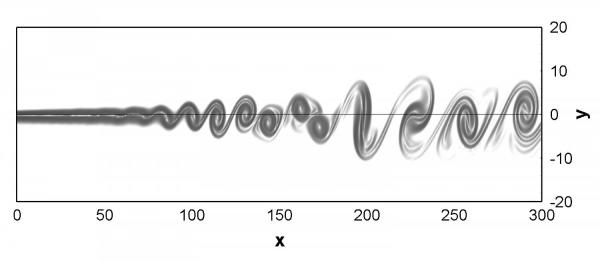}\\
			\includegraphics[width=0.9\textwidth]{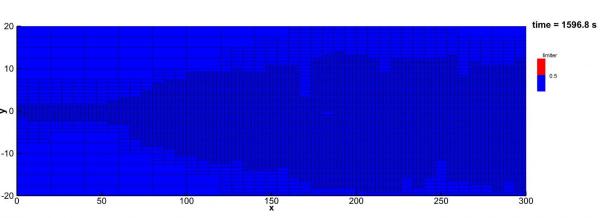}
\end{tabular} 			
			\caption{\textbf{Compressible mixing layer.} Vorticity field (top row) and AMR grid (bottom row) $\omega_z$ obtained with the ADER-DG-$\mathbb{P}_5$ scheme supplemented with the \emph{a posteriori} ADER-WENO SCL	for the compressible 
			mixing layer test for $\mu=10^{-3}$ at $t=68\,T_f=1596.8$s with $T_f=1/f_0$, where $f_0$ is the fundamental frequency of the mixing layer. Up to $\ell_\text{max}=2$ maximum number of refinement  
			levels with  a refine factor $\mathfrak{r}=3$ are used. The limiter is never active. A reference solution \cite{DumbserNSE} for the vorticity field obtained with a high order $P_3P_5$ scheme using 
			a locally refined unstructured triangular grid is provided for comparison (middle row).} \label{fig:CML}
\end{figure}

\begin{figure} 
\centering 
			\includegraphics[width=0.49\textwidth]{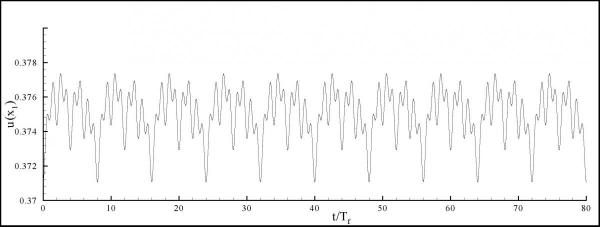}
			\includegraphics[width=0.49\textwidth]{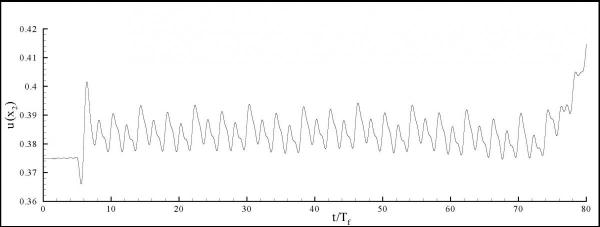}
			\\
			\includegraphics[width=0.49\textwidth]{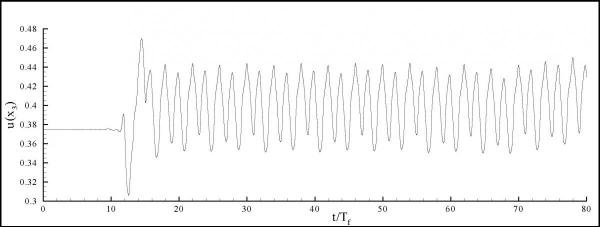}
			\includegraphics[width=0.49\textwidth]{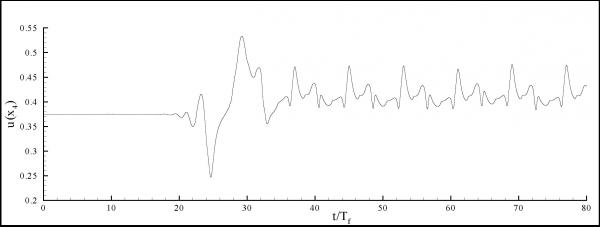}
			\\
			\includegraphics[width=0.49\textwidth]{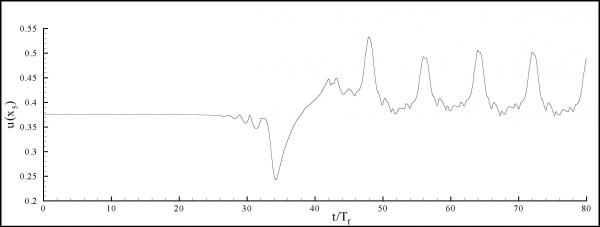}
			\caption{\textbf{Compressible mixing layer.} History of the horizontal velocity component evaluated at $y=0$, along five different axial positions $x_1=0$, $x_2=45$, $x_3=100$, $x_4=200$ and $x_5=285$ (from left to right, from the top to the bottom) for the compressible mixing layer test at  $\mu=10^{-3}$.
			} \label{fig:CML4}
\end{figure}

In this test, originally proposed in \cite{colonius} and then extended to three space dimensions in the work of \cite{babucke} and then reproposed also in \cite{stedg2,DumbserNSE,HPRmodel}, 
the high order of accuracy of our ADER-DG scheme and the judiciousness of the implementation of the SCL are tested. A well known unsteady physical instability is generated along a compressible two  dimensional mixing layer, between the parallel motion of two streams. 
 The upper stream flows at velocity $u_{\infty}=0.5$, the lower one at $u_{-\infty}=0.25$ corresponding to a velocity ratio $\lambda =u_{\infty}/u_{-\infty}=2$; pressure and density are initialized as $\rho=\rho_0=1$ and $p=p_0=1/\gamma$ with a ratio of specific heats $\gamma=1.4$. The singularity at $y=0$ has been smoothed by means of a very simple hyperbolic tangent function 
\begin{align}
u = \frac{1}{8} \left( \tanh (2y) + 3 \right).
\end{align}
 The lengths are made dimensionless with respect to the vorticity thickness at the inflow, given by
\begin{align}
\delta_{\omega_z} (x_0) = \frac{u_{\infty}-u_{-\infty}}{\text{max}( \left. \frac{\partial u}{\partial y}\right|_{x_0})}:= 1,
\end{align}
which allows to define the corresponding Reynolds number
\begin{align}
Re = \frac{\rho_0 u_{\infty} \delta_{\omega_z} (x_0)}{\mu}.
\end{align}
From a rigorous linear stability analysis of the inviscid Rayleigh equations, a proper oscillatory forcing term can be introduced at the inflow in order to facilitate the instability to arise. More details about the resolution of the inviscid Rayleigh but also the viscous Orr-Sommerfeld equations are available in the work of \cite{colonius} and \cite{babucke}.
 Here, the following very simple perturbation has been introduced at the left boundary 
\begin{align}
\delta  (y,t) = A(y) \left[\cos(\omega_0 t) + \cos(\omega_1 t + \phi_1) + \cos(\omega_2 t + \phi_2) + \cos(\omega_3 t + \phi_3)  \right]
\end{align}
where:  $\omega_0=-2\pi f_0 = -0.3147876$ is the fundamental angular frequency ($f_0\approx 0.0501$); $\omega_1=\omega_0/2$, $\omega_2=\omega_0/4$ and $\omega_3=\omega_0/8$ are the corresponding first three subharmonics; $\phi_1=-0.028$, $\phi_2=0.141$ and $\phi_3=0.391$ are the chosen phase-shift of the subharmonic with respect to the fundamental perturbation that allow  to minimize the distance of the vortex pairing, according to \cite{colonius}; $A(y)$ is an amplitude factor ($A<<1$) that can be chosen in the form of a Gaussian distribution centered in the origin of the physical instability, i.e.
\begin{align}
A(y) = \tilde{A} e^{-y^2/4}, \quad \quad \tilde{A} = -10^{-3}.
\end{align} 
The spatial domain is $\Omega = [-50,50]\times[0,400]$, discretized by only $20\times 40$ elements at the coarsest grid level, with a refinement factor $\mathfrak{r}=3$ 
and up to $\ell_\text{max}=2$ maximum number of refinement levels. An ADER-DG-$\mathbb{P}_5$ scheme is employed, supplemented by a third order ADER-WENO finite volume sub-cell limiter. 
Fig.  \ref{fig:CML} 
shows the numerical results for the density variable and the AMR grid by choosing 
a dynamic viscosity of $\mu_1=10^{-3}$
, corresponding to a Reynolds number  of $Re_1=500$.
 The obtained results are directly comparable with the results available in the papers of \cite{colonius,babucke,stedg2,DumbserNSE,HPRmodel} with a good agreement. Fig. 
 \ref{fig:CML} 
 shows the obtained results for the density and vorticity variables, and the AMR grid colored by the limiter-status (limited cells are highlighted in red, unlimited cells are plotted in blue). 
%
 The first vortex pairing occurs at around $x_{p'}\sim 190$. 
 We notice that the SCL has never been activated during the simulation and this is because the physics of the fluid flow has been well-resolved and no spurious oscillations are generated. This is a very  important result and we would like to stress at this point that the presented sub-cell limiting procedure does \textit{not} dissipate the real physical instabilities, but only the numerical ones,  preserving the original resolution of high order unlimited DG scheme for smooth flows. 
Finally, Fig.  \ref{fig:CML4} shows the comparison of the time series of the horizontal velocity evaluated at $y=0$ at different axial positions. These plots give a better idea on the time-scales of the  development of the instability; they seem to be well compatible with literature results  (see \cite{colonius,DumbserNSE}).

\subsection{Shock-vortex interaction}
\label{test:SVI}
\begin{figure} 
\centering 
			\includegraphics[width=\textwidth]{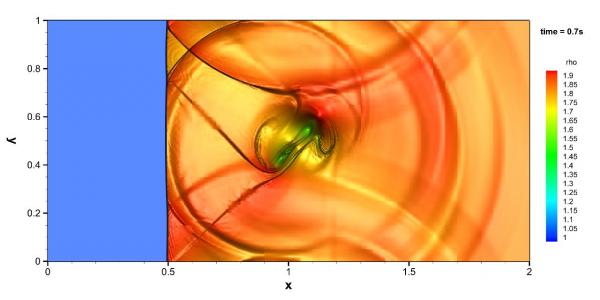}\\ 
			\includegraphics[width=\textwidth]{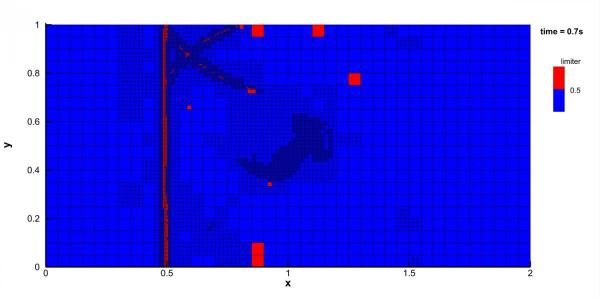} 
			\caption{\textbf{Shock-vortex interaction.} Density (top) and AMR grid colored by the limiter status (bottom) obtained with our ADER-DG-$\mathbb{P}_5$ supplemented with the \emph{a posteriori} TVD SCL in \emph{primitive variables}	for the shock-vortex interaction test at $t=0.7$s. Up to $\ell_\text{max}=2$ maximum number of refinement levels with  a refine factor $\mathfrak{r}=3$ are used. The kinematic viscosity is $\nu = 10^{-8}$.} \label{fig:SVI1}
\end{figure}

\begin{figure} 
\centering 
			\includegraphics[width=\textwidth]{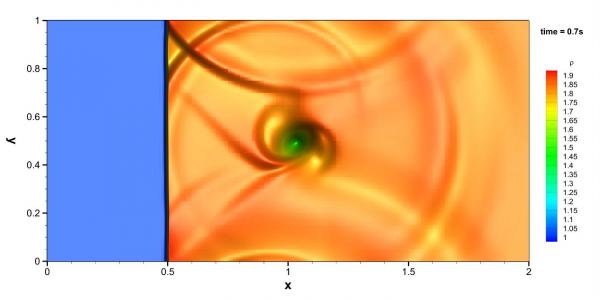}\\ 
			\includegraphics[width=\textwidth]{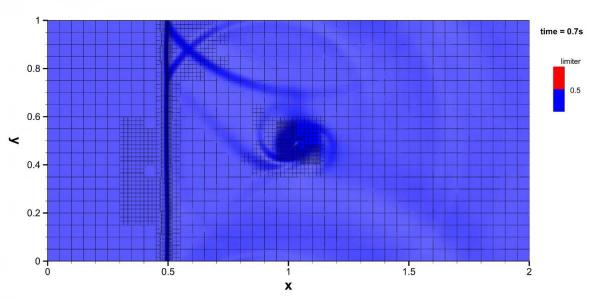}
			\caption{\textbf{Shock-vortex interaction.} Density (top),  and limiter status (bottom) obtained with our ADER-DG-$\mathbb{P}_5$ supplemented with the \emph{a posteriori} SCL	for the shock-vortex interaction test at $t=0.7$s. Up to $\ell_\text{max}=2$ maximum number of refinement levels with  a refine factor $\mathfrak{r}=3$ are used. The kinematic viscosity is $\nu = 10^{-3}$.} \label{fig:SVI2}
\end{figure}

An interesting two dimensional problem for testing the AMR framework dealing with shocks and smooth waves together is the so called \emph{shock-vortex interaction} test. In this problem a smooth vortex hits a stationary normal shock wave, representing an optimal scenario for testing high order shock capturing schemes. The spatial domain is $\Omega=[0,2]\times[0,1]$ with periodic boundary conditions in the vertical direction, analytical boundary conditions at the left boundary and a classical outflow boundary condition at the right. The vortex is centered at $(x_V,y_V) = (0.25,0.5)$ and its strength is characterized by a Mach number of $M_V= v_m/c_0 =0.7$, $c_0=\sqrt{\gamma p_0/\rho_0}$ being the adiabatic sound speed upstream the shock, with $p_0=1$ and $\rho_0=1$. The angular velocity $\omega_V$ is distributed according to
\begin{equation}\label{eq:SVI}
  \omega_V = \left\{\begin{array}{ll}
\omega_m \frac{r}{a} & \;\textrm{for}\quad r \leq a\,, \\
 \noalign{\medskip}
 \omega_m \frac{a}{a^2-b^2}\left(r-\frac{b^2}{r}\right) & \;\textrm{for}\quad a\leq  r \leq b\,, \\
 \noalign{\medskip}
0 & {\rm otherwise} \,,\\
 \end{array}\right. 
\end{equation} 
where $r^2=(x-x_V)^2 + (y-y_V)^2$. Pressure and density are evaluated according to the equations
\begin{eqnarray}
p=p_0\left(\frac{T}{T_0}\right)^{\frac{\gamma}{\gamma-1}}, \quad \quad
\rho=\rho_0\left(\frac{T}{T_0}\right)^{\frac{1}{\gamma-1}}\,.
\end{eqnarray}
after solving the ordinary differential equation for the temperature
\begin{equation}
\frac{dT}{dr}=\frac{\gamma-1}{R\gamma}\frac{\omega_V^2(r)}{r}\,.
\end{equation}
The  unperturbed upstream variables are chosen in compliance with the equation of state of ideal gases $p_0=R\rho_0T_0$ where the the gas constant is $R=1$. 
The remaining parameters are chosen to be $\gamma=1.4$, $a=0.0075$ and $b=0.175$ and the Prandtl number of $Pr=0.7$.
Finally, the stationary shock with Mach number $M_S=1.5$ is placed at $x=0.5$ and the downstream variables are computed according to the classical Rankine-Hugoniot conditions \cite{Landau-Lifshitz6}.
The current test has been solved with the $\mathbb{P}_5$ version of our ADER-DG method, supplemented only by a second order accurate shock capturing TVD finite volume scheme on the subgrid, 
based on reconstruction in primitive variables. Fig.  \ref{fig:SVI1} shows the computed results for the density variable and the AMR grid colored by the limiter status for $\mu=10^{-8}$ (limited cells are highlighted in red, unlimited cells are plotted in blue). The obtained results are in agreement with the results available in literature \cite{Dumbser2014,Dumbser2007,donat}. Moreover, Fig.  \ref{fig:SVI2} shows the  computed results obtained by choosing a viscosity of $\mu=10^{-3}$. The effect of higher physical viscosity is evident, since the final solution is much smoother because of the presence of viscous effects and
heat conduction. We observe that the SCL is never activated (although the entire MOOD framework is switched on in this test problem!), 
since for sufficiently resolved viscous flows, the use of a limiter becomes unnecessary. 

\subsection{Viscous double Mach reflection problem}
\label{test:vDMR}
\begin{figure} 
\centering 
			\includegraphics[width=0.4\textwidth]{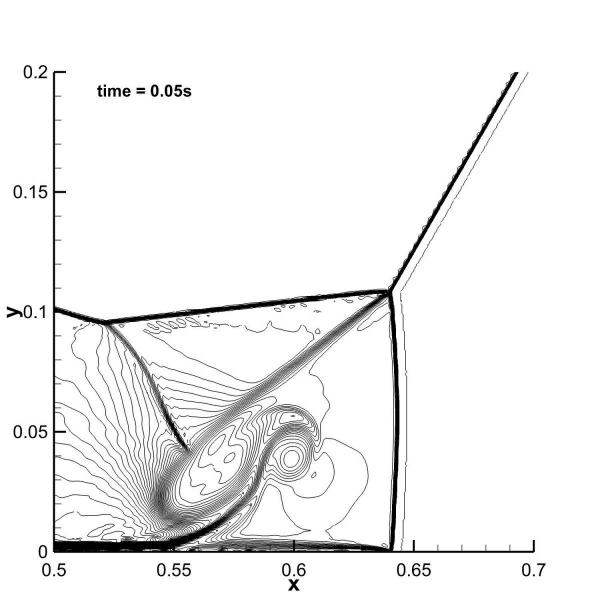}  
			\includegraphics[width=0.4\textwidth]{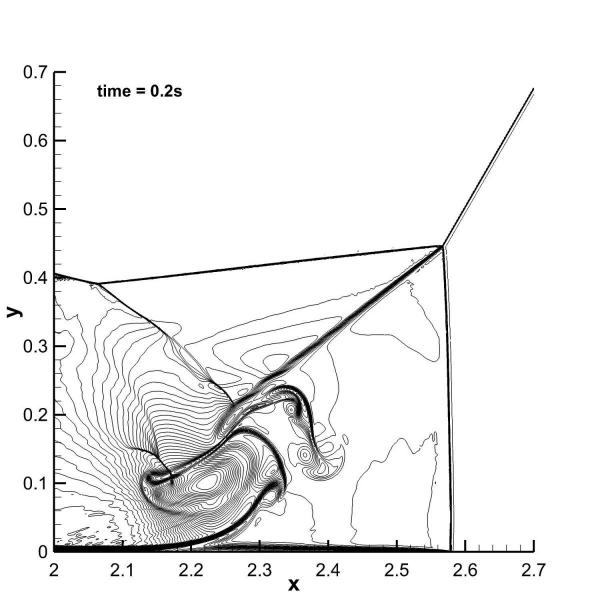} \\
			\includegraphics[width=0.4\textwidth]{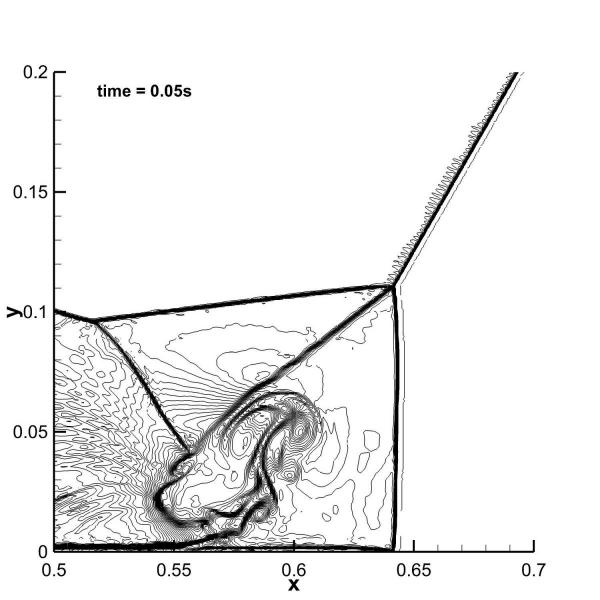}  
			\includegraphics[width=0.4\textwidth]{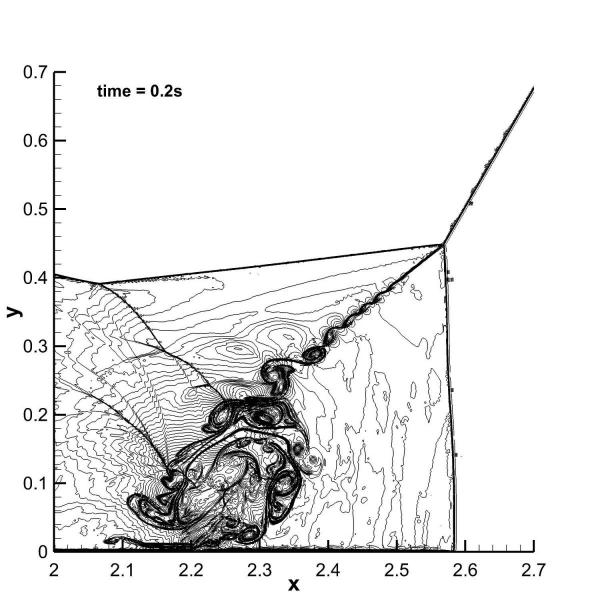}  \\
			\includegraphics[width=0.4\textwidth]{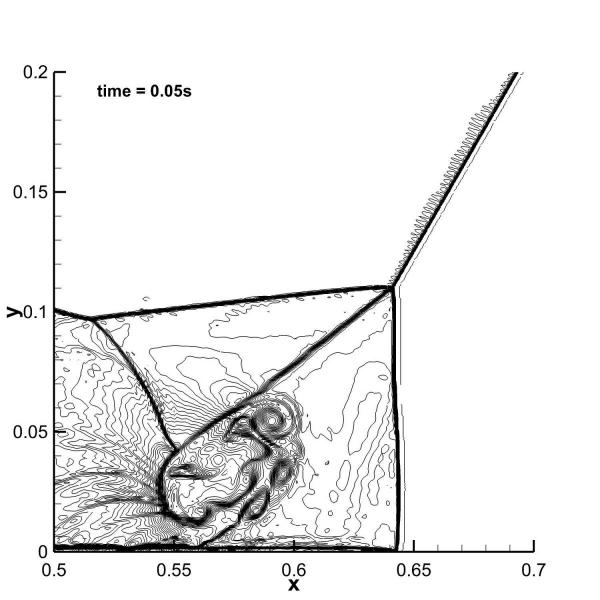}  
			\includegraphics[width=0.4\textwidth]{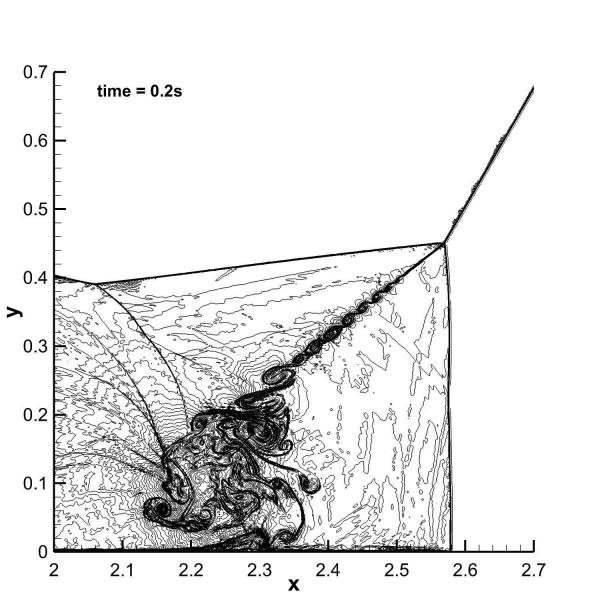}   
			\caption{\textbf{Viscous double Mach reflection.} Contour lines of the density for the viscous double Mach reflection test for viscosity $\mu=10^{-3}$, $10^{-4}$ and  the inviscid limit  $10^{-8}$, from top to bottom, at different times 
			$t=0.05$ (left) and $0.20$ (right), obtained with ADER-DG-$\mathbb{P}_5$ and \emph{a posteriori} SCL WENO$3$.} \label{fig:vDMR_lines_1}
\end{figure}

Originally proposed by Woodward and Colella in \cite{Woodward84} for the inviscid case, here we solve a \textit{viscous} version of the two dimensional double Mach reflection problem at very high 
Mach number ($M_S=10$). In this test, a planar shock wave hits a rigid wall at an angle of incidence of $\alpha_S=60^{\circ}$. The shock wave reflection, the viscous wall boundary layer, but also the  physical instabilities inside the front of incidence for high Reynolds numbers, make this scenario very intriguing for testing the ability of a high-order numerical scheme to capture all the flow 
physics, from the smaller and the larger vortex structures to the strong shock waves that appear at $M=10$, (see \cite{Woodward84,DumbserNSE,Dumbser2014,Zanotti2015c}). 
The $60^{\circ}$-inclined wavefront of the viscous shock-wave is initialized by imposing $x=0$ as the initial point of incidence at the wall, 
and prescribing the classical Rankine-Hugoniot conditions of the compressible Euler equations at the shock interface with respect to the chosen downstream variables, having 
\begin{align}
\left( \rho_0, u'_0,v'_0,p_0 \right) = \left\{ \begin{array}{ll}  \left( 8 , 8.25 \gamma, 0 , 116.5 \right)\frac{1}{\gamma} &  \text{for} \,\, x'\,\, \text{upstream} \\ \left( \gamma , 0, 0 , 1 \right)\frac{1}{\gamma} &  \text{for} \,\, x'\,\, \text{downstream} \end{array} \right. . 
\end{align}
where the primed variables and coordinates $u'$, $v'$ and $x'$ are evaluated with respect to the rotated coordinate system, $x'$ being the streamwise direction.  
 The Prandtl number is $Pr=3/4$. The spatial domain is chosen  to be $\Omega=[0,4]\times [0,1]$ with no-slip boundary condition at the bottom, outflow boundary condition at the right, and the aforementioned analytical solution of the moving incident shock-wave in the remaining left and top boundaries. 
For this test the ADER-DG-$\mathbb{P}_5$ scheme has been used, together with a third order ADER-WENO finite volume scheme as subcell limiter (SCL). 
The coarsest mesh, the one of the $0-th$ refinement level, is made up of $80\times 20$ elements upgraded by up to $\ell_\text{max}=2$ maximum number of refinement levels with  a refine factor $\mathfrak{r}=3$. Then, the corresponding characteristic lengths $h_\ell$ of the three refinement levels are $h_0=1/20$, $h_1=1/60$ and $h_2=1/180$ and the \emph{effective characteristic lengths} $\tilde{h}_\ell$ that take account of the $d.o.f$ of the polynomial basis ($\tilde{h}_\ell = h_{\ell}/(N+1)$) are $\tilde{h}_0 = 1/120 $, $\tilde{h}_1 = 1/360$ and  $\tilde{h}_2 = 1/1080$.
Figure \ref{fig:vDMR_lines_1}   shows the numerical results for the density contour lines at time $t=0.05$ and $0.2$ for differed dynamic viscosity coefficients, i.e. $\mu_1=10^{-3}$ leading to the shock-Reynolds number $Re_1=\rho_0 M_S /\mu_1=10^4$, $\mu_2=10^{-4}$ leading to $Re_2=10^5$  and the almost inviscid limit case $\mu_3 =10^{-8}$. 
It is important to note that when the inviscid compressible Euler equations are solved, the present problem will develop smaller and smaller spatial scales in an unbounded manner, since there is no physical
viscosity in the Euler equations that prevents the generation of small scale vortex structures. In the presence of physical viscosity, however, there exists a smallest spatial scale at which vortex structures dissipate energy into internal energy and below which no smaller spatial scales can exist. 
%

\begin{figure} 
\centering 
			\includegraphics[width=0.4\textwidth]{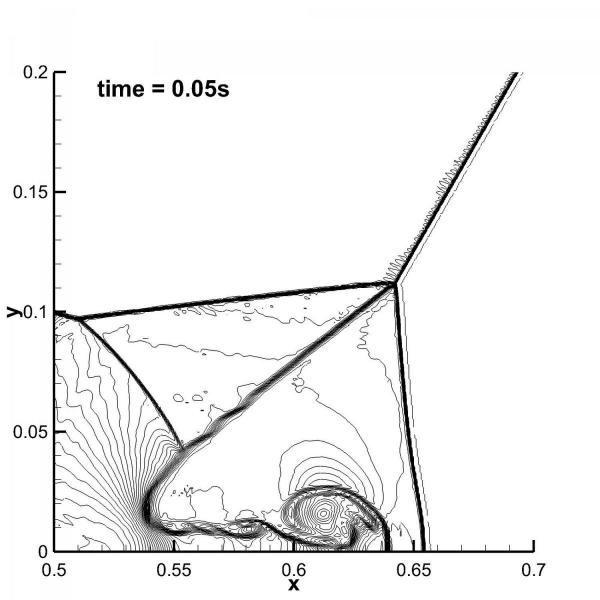}  
			\includegraphics[width=0.4\textwidth]{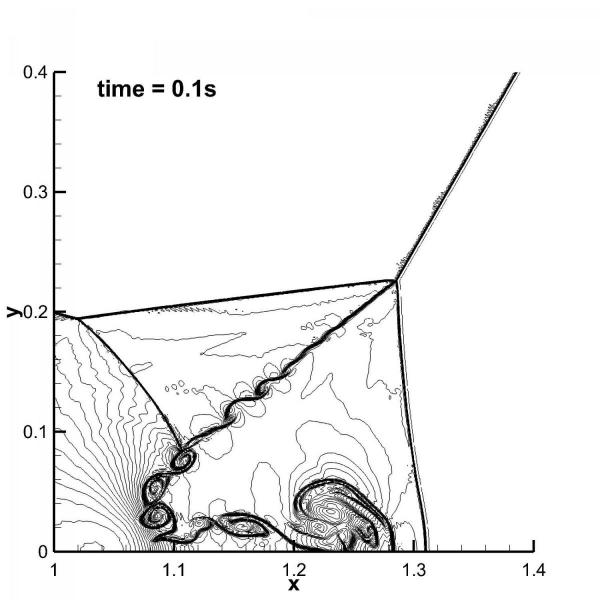} \\
			\includegraphics[width=0.4\textwidth]{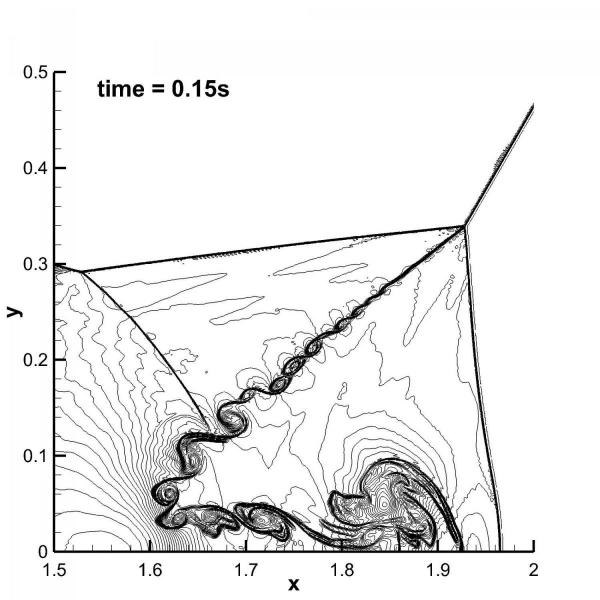}  
			\includegraphics[width=0.4\textwidth]{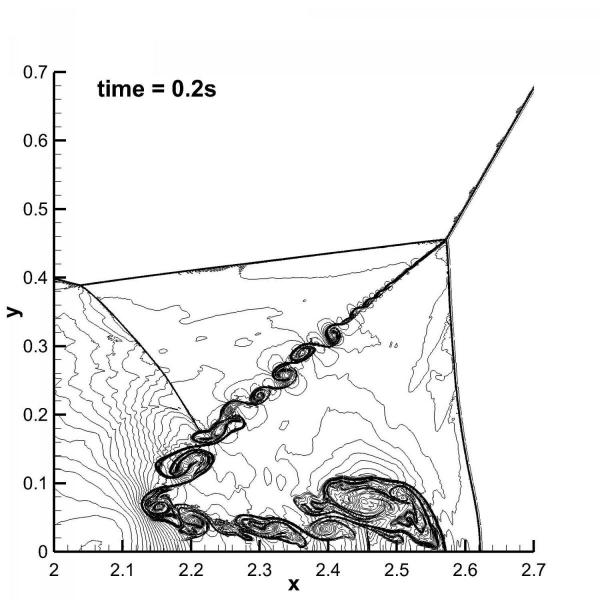}   
			\caption{\textbf{Viscous double Mach reflection.} Contour lines of the density for the viscous double Mach reflection test \textbf{with purely reflective wall boundary conditions} for viscosity $10^{-4}$ at different times 
			$t=0.05$, $0.10$, $0.15$ and  $0.20$,  from top left to bottom right, obtained with ADER-DG-$\mathbb{P}_5$ and \emph{a posteriori} SCL WENO$3$.} \label{fig:vDMR_lines_3}
\end{figure}

\begin{figure} 
\centering 
			\includegraphics[width=0.95\textwidth]{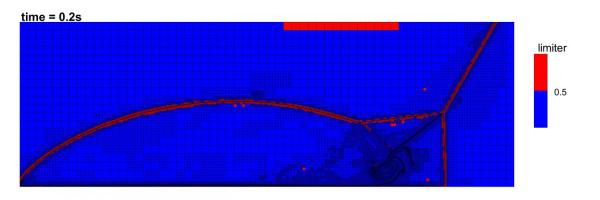}  \\
			\includegraphics[width=0.95\textwidth]{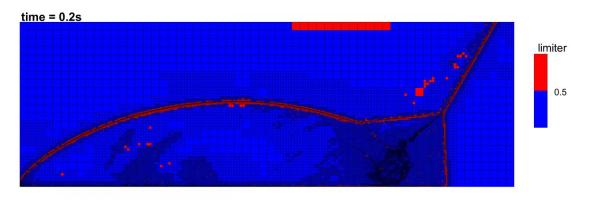}  \\
			\includegraphics[width=0.95\textwidth]{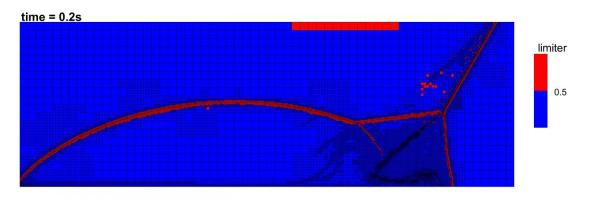} \\ 
			\caption{\textbf{Viscous double Mach reflection.} Plot of the AMR grid for ADER-DG-$\mathbb{P}_5$ polynomials (blue) and the ADER-WENO$3$ \emph{sub-cell} averages, i.e. the limited cells (red),  for the viscous double Mach reflection 
			test at the final time $t=0.2$ obtained by choosing  $\mu=10^{-3}$ and no-slip walls (top), $\mu=10^{-4}$ and no-slip walls (center) and $\mu=10^{-4}$ with slip walls (bottom).  
			} \label{fig:vDMR_limiter}
\end{figure}

The classical '\emph{crow's \reds{feet}}'-shaped  (i.e. the right $3+1$ shock-wave-interfaces that are \emph{incident} with respect to a central node) 
front-wave is well reproduced.  It holds some interest noticing the differences of the distance 
 between the central node and the location of the first vortex appearance along the central slip line, which is affected by a Kelvin-Helmholtz instability at higher Reynolds numbers. 
%
%
Then, the same numerical simulation has been repeated by applying reflective (inviscid) slip-wall boundary conditions at the bottom, instead of the classical no-slip wall-boundary conditions.
Fig.  \ref{fig:vDMR_lines_3} shows the numerical solution for the density contour lines obtained  
at time $t=0.2$. 
 Notice that the no-slip boundary conditions at the bottom wall lead to a completely different flow pattern compared to the usual slip wall boundaries used for the simulation of inviscid flows: 
the development of the well known 'mushroom'-type shape of the the purely reflective slip-wall case  is prevented  because of the thin boundary layer at the wall, leading to $\partial u/\partial y \neq 0$ at $y=0$. 
The complete AMR grids colored by the limiter status are depicted in Fig. \ref{fig:vDMR_limiter} for the considered Reynolds number regimes and boundary conditions. One can notice that the AMR 
method worked properly, following the main shock waves and resolving also the vortexes generated by the Kelvin-Helmholtz instability along the slip line. Moreover, also the SCL ADER-WENO$3$ is essentially activated only when and  where it is necessary, i.e. only next to the stronger shocks (see red cells in Fig. \ref{fig:vDMR_limiter} allowing the ADER-DG $\mathbb{P}_{5}$-polynomials to represent the numerical solution in the smoother zones and throughout the non-linear instabilities. Notice that only a minor number of '\emph{false-positive}' limited cells have been detected for this test-problem. 

It should be emphasized that there are not many reference results published in the literature concerning the \textit{viscous} double Mach reflection problem. In the case of high 
Reynolds numbers and inviscid slip wall boundary conditions, our obtained results seem to be in good agreement with the results present in the literature  
\cite{Woodward84,DumbserNSE,Dumbser2014,Zanotti2015c,HPRmodel}. 

\subsection{Kelvin Helmholtz instability for the CNS and the VRMHD equations}
\label{test:KH}
\begin{figure}
\centering 
			 \includegraphics[width=0.8\textwidth]{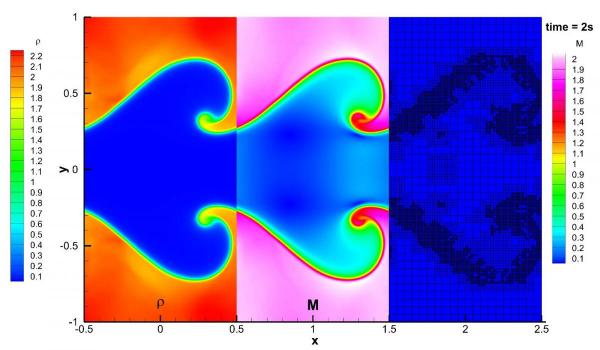} 
			\\ \includegraphics[width=0.8\textwidth]{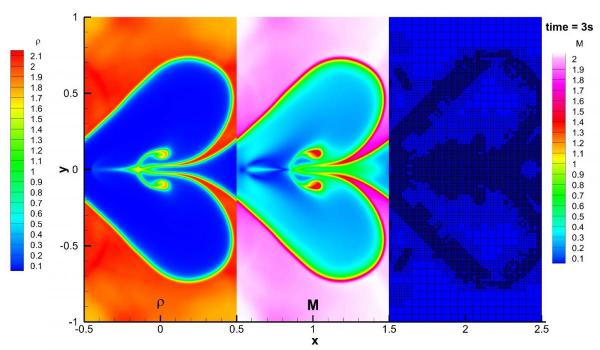} 
			\\ \includegraphics[width=0.8\textwidth]{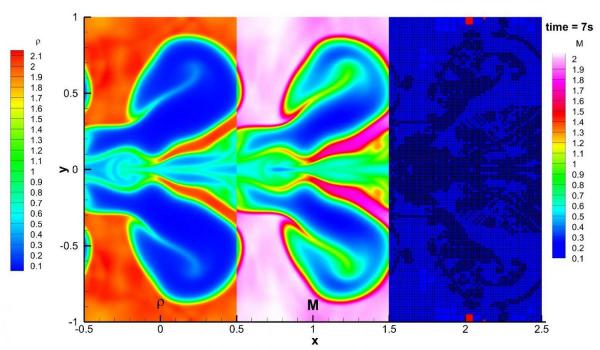} 
			\caption{\textbf{Kelvin Helmholtz instability (Navier-Stokes).} Numerical solution of the compressible Navier-Stokes equations for the two dimensional Kelvin-Helmholtz instability our ADER-DG-$\mathbb{P}_{3}$ supplemented by the \emph{a posteriori} SCL using $20\times 40$ elements on the coarsest level, up to $\ell_\text{max}=2$ maximum number of refinement levels with  a refine factor $\mathfrak{r}=3$. The density (left), the local Mach number (center) and the   active-mesh colored by the limiter-status are plotted at times $t=2.0$, $3.0$ and $7.0$ from the top to the bottom, respectively.}\label{fig:KHNS}
\end{figure}

\begin{figure} 
\centering 
			 \includegraphics[width=0.8\textwidth]{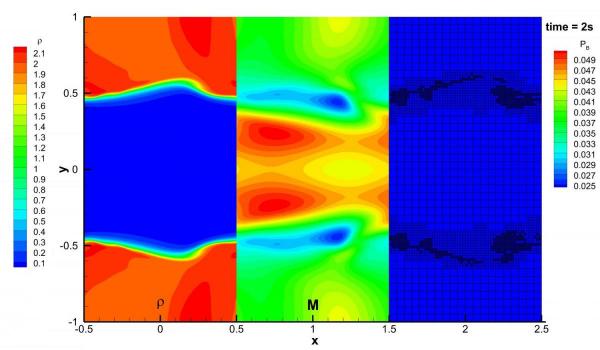} 
			\\ \includegraphics[width=0.8\textwidth]{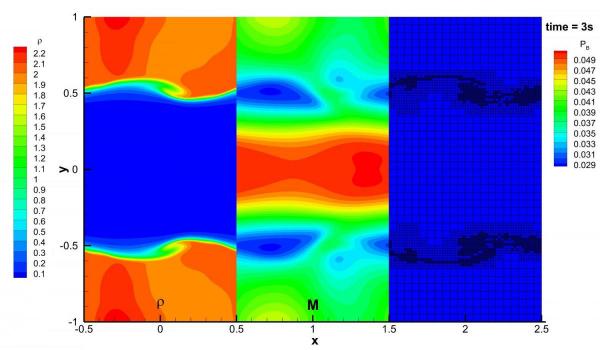} 
			\\ \includegraphics[width=0.8\textwidth]{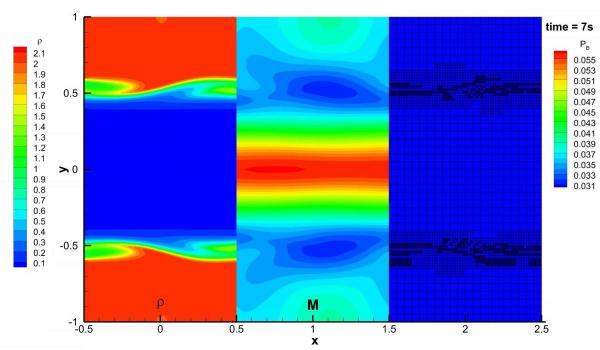} 
			\caption{\textbf{Kelvin Helmholtz instability (resistive MHD).} Numerical solution of the viscous and resistive MHD equations for the two dimensional Kelvin-Helmholtz instability our ADER-DG-$\mathbb{P}_{3}$ supplemented by the \emph{a posteriori} SCL using $20\times 40$ elements on the coarsest level, up to $\ell_\text{max}=2$ maximum number of refinement levels with  a refine factor $\mathfrak{r}=3$.
			The density (left), the magnetic pressure $|\mathbf{B}|/8\pi$ (center) and the active-mesh colored by the limiter-status are plotted at 
			times $t=2.0$, $3.0$ and $7.0$ from the top to the bottom, respectively.}\label{fig:KHMHD}
\end{figure}

In this two-dimensional test the well known physical instability that takes the name from William Thomson (named Lord Kelvin) and Hermann von Helmholtz is simulated both for the compressible Navier-Stokes and the viscous-resistive MHD equations. The Kelvin-Helmholtz instability plays important roles in dissipative processes and momentum/energy transfer in atmospheric processes, fluvial engineering, oceanography, but also solar physics and astrophysics. In general, it is the physical instability that arises in the nonlinear interaction of the relative motion of two parallel fluids, as in the compressible mixing layer problem solved before. The spatial domain $\Omega=[-0.5,0.5]\times[-1,1]$ is discretized on the zeroth level with only $20\times 40$ elements. The AMR framework is used 
up to $\ell_\text{max}=2$ maximum number of refinement levels and a refine factor $\mathfrak{r}=3$. Periodic boundary conditions are assumed at the borders. The fluid flow is initialized following \cite{Mignone2009}, \cite{Beckwith2011}, \cite{Radice2012a} and \cite{Zanotti2015d}, i.e.
\begin{equation}\label{KHI-vx}
  u = \left\{\begin{array}{ll}
  v_s \tanh{[(y-0.5)/a]} & \quad y > 0\,, \\
 \noalign{\medskip}
 -v_s \tanh{[(y+0.5)/a]}  & \quad y \leq 0    \,, \\
 \noalign{\medskip}
 \end{array}\right.
\end{equation}
where $v_s=1.0$ is the velocity of the shear layer and $a=0.01$ is its characteristic size. A small
transverse velocity has been conveniently introduced to trigger the instability by choosing
\begin{equation}\label{KHI-vy}
  v = \left\{\begin{array}{ll}
  \eta_0 v_s \sin{(2\pi x)} \exp{[-(y-0.5)^2/\sigma]} & \quad y > 0\,, \\
 \noalign{\medskip}
 -\eta_0 v_s \sin{(2\pi x)} \exp{[-(y+0.5)^2/\sigma]}  & \quad y \leq 0    \,, \\
 \noalign{\medskip}
 \end{array}\right.
\end{equation}
with $\eta_0=0.1$ and $\sigma=0.1$. Finally, the fluid density is
\begin{equation}\label{KHI-rho}
  \rho = \left\{\begin{array}{ll}
  \rho_0 + \rho_1 \tanh{[(y-0.5)/a]} & \quad y > 0\,, \\
 \noalign{\medskip}
 \rho_0 - \rho_1 \tanh{[(y+0.5)/a]}  & \quad y \leq 0    \,, \\
 \noalign{\medskip}
 \end{array}\right.
\end{equation}
with $\rho_0=1.005$ and $\rho_1=0.995$. The dynamic viscosity coefficient has been chosen to be $\mu=10^{-3}$. For the MHD case the electric resistivity is $\eta=10^{-2}$ and a constant magnetic field is initialized horizontally oriented as
\begin{align}
\left( B_x, B_y, B_z \right) = \left( 0.1, 0, 0 \right).
\end{align}
Figures \ref{fig:KHNS} and \ref{fig:KHMHD} show the numerical results obtained with our ADER-DG-$\mathbb{P}_{3}$ scheme supplemented with the \emph{a posteriori} sub-cell WENO$3$ limiter for the compressible NS and the resistive MHD equations, respectively, up to the time $t_e=7$. It becomes evident how the initial magnetic field drastically influences the dynamics of the electrically conducting fluid (see Fig. \ref{fig:KHMHD}). At first, the hydrodynamical forces between the two layers are in mutual unstable equilibrium. By introducing the nonzero vertical velocity component we upset the balance and the fluid state starts falling, looking for a new equilibrium state through the generation of mixing-breaking waves and diffusion processes (see Fig. \ref{fig:KHNS}). When a non-negligible magnetic field is active, every minimal distortion in the fluid flow corresponds to a deformation in the magnetic field. In this sense an amount of work is necessary to the magnetic-field lines to distort, and therefore also to the streamlines. Consequently, the resulting mixing process is weakened with respect to the non-charged fluid flow.  Notice how a non-negligible magnetic pressure gradient pushes the  fluid flow from the inner lower-density core to the outer zones (see Fig. \ref{fig:KHMHD}), causing the shear layer to remain spatially confined for longer times. 
Thus the magnetic field plays the role of stabilizer of the initial unstable equilibrium, leading to a longer life-time of the double shear layer flow.

\subsection{Magnetic reconnection}
\label{mag-rec}
In this test we consider the classical problem of magnetic reconnection,
which consists of the re-adjustment of the magnetic field topology
due to a non-vanishing resistivity,
typically occurring through sheet-like structures of length $L$ and width $a$.
The classical Sweet--Parker (SP)
reconnection model predicts a dissipation of magnetic energy with
a reconnection timescale $\tau_{\rm rec}\sim\tau_A S^{1/2}$, where $S$ is the Lundquist number.
Since both in astrophysical context and in laboratory conditions the Lundquist number
is very large ($S\sim 10^{12}$ in the solar corona and $S\sim 10^{8}$ in tokamaks), the interest towards simple resistive MHD reconnection has been 
frustrated for a long time. However, a novel attention has been triggered by the discovery that
current sheets with large aspect ratios $L/a$ become violently unstable \cite{Biskamp1986,Loureiro2007,Samtaney2009, Landi2015}, 
generating plasmoid chains on smaller and
smaller scales. 

\begin{figure}
\centering 
			 \includegraphics[width=\textwidth]{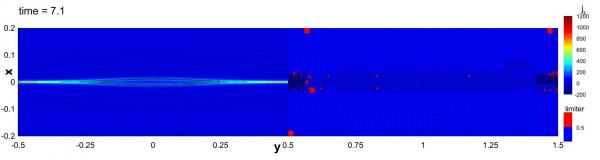} \\
			 \includegraphics[width=\textwidth]{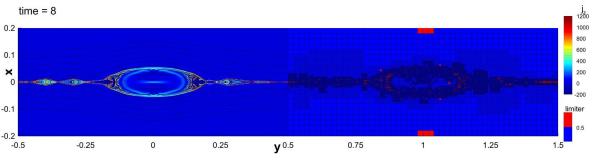} \\
			 \includegraphics[width=\textwidth]{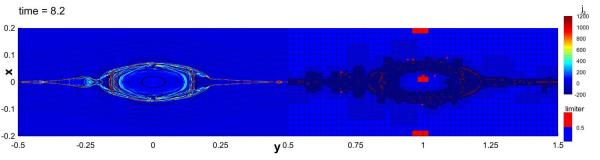} \\
			 \includegraphics[width=\textwidth]{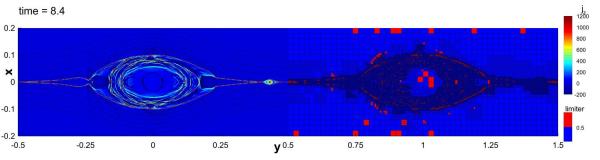} \\
			 \includegraphics[width=\textwidth]{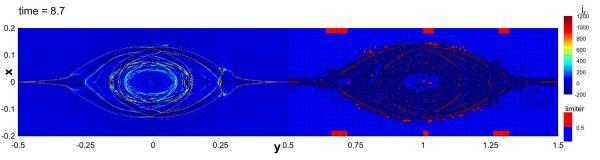} 
			\caption{\textbf{Magnetic reconnection.} Numerical solution of the resistive MHD equations for the two dimensional magnetic reconnection test problem at several time-step  obtained with our ADER-DG-$\mathbb{P}_{5}$ supplemented by the \emph{a posteriori} WENO$3$ SCL using $20\times 50$ elements on the coarsest level, up to $\ell_\text{max}=2$ maximum number of refinement levels with  a refine factor $\mathfrak{r}=3$. The density current $j_z$ (left) and the   active-mesh colored by the limiter-status (right) are plotted at times $t=7.1$,  $8.0$, $8.2$,  $8.4$ and $8.7$ from the top to the bottom, respectively.}\label{fig:mRec}
\end{figure}

\begin{figure}
\centering 
			 \includegraphics[width=\textwidth]{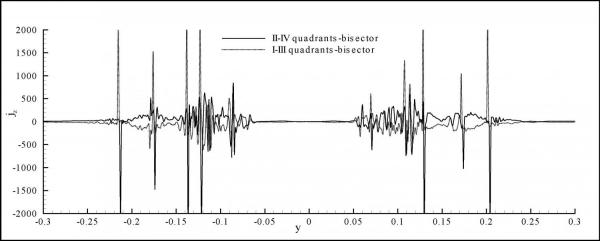}\\
			 \includegraphics[width=\textwidth]{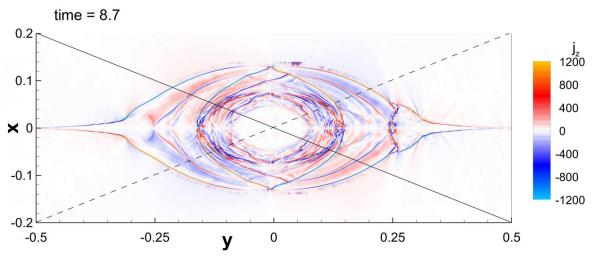}
			\caption{\textbf{Magnetic reconnection.} Interpolation of the current density $j_z$ for the two dimensional magnetic reconnection test problem at time $t=8.7$ along the two bisectors of the rectangular computational domain (top), highlighting  the \emph{tree-ring} structure of the alternately-positive$/$negative current density of the major plasmoid. At the bottom the corresponding current density $j_z$ is shown and the two considered bisector have been highlighted in white continuous and dashed lines.}\label{fig:Jz_bisectors}
\end{figure}

Here we reproduce a representative case of magnetic reconnection with our ADER-DG scheme, focusing on the {\em ideal tearing mode} investigated recently by \cite{Landi2015}. 
The numerical domain is  $[-20a,20a]\times[-L/2,L/2]$, where $a=L/S^{1/3}$ is the width of the current sheet, while the
Lundquist number $S$, which is given by the ratio between the diffusion timescale $\tau_D=L^2/\eta$ and 
the advection timescale $\tau_A=L/v_a$, is $S=L v_a/\eta$. The magnetic field in the $(x,y)$ plane follows the typical Harris model, with, in addition, 
a perpendicular component, in order to have a globally uniform magnetic field at time $t=0$, i.e.
\begin{equation}
{\bf B}=B_0\left[ \tanh(x/a) {\bf \hat y} + \sech(x/a)  {\bf \hat z} \right]\,,
\end{equation}
where $B_0$ is related to the Alfven-speed by the usual expression $v_a^2=B_0^2/(4\pi \rho)$.
The thermal pressure, which is also initially uniform over the computational domain,
is determined through a condition on the magnetic Mach number $M=v_a/c_s$. For an ideal gas equation of state $p=\rho\epsilon(\gamma-1)$, 
this allows to obtain $p=\rho/(\gamma M^2)$.
In our test we have chosen $v_a=L=1$, $\gamma=5/3$, $M=0.7$ and $S=10^6$, corresponding to 
a current sheet thickness $a=0.01$ and to
an asymptotic  plasma parameter $\beta=2.4$.
Like in \cite{Landi2015}, the instability is triggered by inserting a perturbation in the
velocity field at time $t=0$, i.e.
\begin{eqnarray}
v_x&=&\varepsilon \tanh\xi \exp(-\xi^2) \cos(ky) \\
v_y&=&\varepsilon (2\xi \tanh\xi - \sech^2\xi)\exp(-\xi^2)S^{1/2}\sin(ky)/k\,,
\end{eqnarray}
where $\varepsilon=10^{-3}$, $\xi=x S^{1/2}$, while the wave-number is computed from $k L=2\pi m$, with $m=10$.
Free outflow and periodic  boundary conditions are chosen along $x$ and $y$, respectively.
The time evolution of the numerical solution for the density current $j_z = \partial_x B_y - \partial_y B_x$ obtained with our ADER-DG-$\mathbb{P}_{5}$ supplemented by the \emph{a posteriori} WENO$3$ SCL is plotted in Fig. \ref{fig:mRec} next to the active-mesh contour plot. The computational domain has been discretized between $20\times 50$ coarsest elements, up to $\ell_\text{max}=2$ maximum number of refinement levels with a refine factor $\mathfrak{r}=3$.
The initial condition consists in a  positive (exiting)  current density $j_z$ localized within a thin vertical  layer centered in $x=0$. Because of Ampere's law, there is a magnetic tension acting along the thin current density (along y) and, therefore, this system can be seen as a \emph{tighten string} that owes its instability to the compressible nature of the fluid. The present test is often referred indeed to as the 'tearing instability' process. 
Due to the initial perturbation, the symmetry of the system breaks and a higher current density-segment follows up next to a lower one. Simultaneously, the magnetic field aims to maintain the divergence free condition and the lower current density-segment is consequently \emph{bifurcated} (see the first plot in Fig. \ref{fig:mRec} keeping in mind the periodic boundary conditions). In this way, the first main \emph{reconnection island} (or \emph{major plasmoid}) is generated and it takes the form of an \emph{harmonic perturbation} of the current density $j_z$.
Then, the higher current density-segment behaves like a source of new smaller reconnection islands that are attracted to the center of the major plasmoid. Throughout this non-linear process the successively generated smaller islands  collide and merge with the major plasmoid, leading to the so called \emph{plasmoid coalescence}. The major plasmoid broadens out, resulting in a larger \emph{onion} like structure of alternating positive$/$negative current density interfaces (see Fig. \ref{fig:mRec} and \ref{fig:Jz_bisectors}).

\section{Conclusions}
\label{sec.concl}

Their high order of accuracy combined with their locality (no reconstruction step needed) make DG methods very attractive for solving a wide range of spatial flow scales in fluid dynamics 
when used within an adaptive mesh refinement framework (AMR). 
However, it is a well known fact that pure DG methods are unable to properly resolve discontinuous waves or very sharp flow profiles without introducing unphysical spurious oscillations 
(Gibbs phenomenon).  
To cope with this problem, in this paper an arbitrary high-order unlimited ADER-DG method has been supplemented with a high-order accurate and very robust ADER-WENO finite-volume method. 
The adopted limiting method, based on an \emph{a posteriori} survey of troubled zones and, when necessary, a complete re-computation of the solution by means of a more robust finite volume 
scheme within a proper finer sub-grid, i.e. the SCL, has been introduced for the first time by \cite{Dumbser2014}. 
The primordial version of the adopted \textit{a posteriori} limiting method is due to the series of papers concerning the multi-dimensional optimal order detection (MOOD) criteria for 
finite-volume methods introduced in \cite{CDL1,CDL2,CDL3,ADER_MOOD_14}. 
The SCL procedure for DG methods has been extended to auto-adaptive meshes by \cite{Zanotti2015c,Zanotti2015d}, but only for inviscid fluids. 
In this work, the cited numerical method has been extended for the first time to solve the fluid dynamics of dissipative flows, i.e. the \emph{viscous} compressible Navier-Stokes equations 
and the \emph{resistive} magneto-hydrodynamic equations (VRMHD). 

The numerical method has been thoroughly tested on a large set of non-trivial numerical benchmark problems, from low to high Mach number flows, from low to high Reynolds number regimes, 
for which a reference solution or published reference results exist. In particular, the higher order of accuracy combined with the shock-capturing capabilities of the method have been 
successfully demonstrated. 

Future research will concern the extension of the present numerical formulation to the equations of \emph{resistive relativistic MHD}, where a special treatment \emph{stiff source terms} 
becomes necessary. 

Moreover, in order to cope with the severe time restriction of \emph{explicit} DG methods, a semi-implicit time discretization seems to be advantageous. 
Very recently, a novel semi-implicit DG approach on \textit{staggered} grids has been introduced, first for the shallow water equations in \cite{DumbserCasulli2013,TavelliDumbser2014}, and 
subsequently reformulated also for the incompressible Navier-Stokes equations in two and three-space dimensions for unstructured (see \cite{TavelliDumbser2014,TavelliDumbser2014b,TavelliDumbser2015}) 
 Cartesian meshes (see \cite{FambriDumbser}), \reds{and also with adaptive-mesh refinement (see \cite{AMRDGSI})}, assuring spectral convergence within a high order space-time DG framework. 
The semi-implicit time-discretization allows the use of large time-steps, while the staggered grid leads to well-conditioned symmetric positive definite linear algebraic systems with small 
computational stencils and smaller linear systems to be solved compared to the same method used on collocated grids. 
Therefore, further work will also concern the incorporation of \textit{staggered} semi-implicit DG schemes into the high order AMR framework used in the present paper. 
 
\section*{Acknowledgments}
The presented research has been financed by the European Research Council (ERC) under the European Union's Seventh Framework 
Programme (FP7/2007-2013) with the research project \textit{STiMulUs}, ERC Grant agreement no. 278267. 

The authors have also received funding from the European Union's Horizon 2020 Research and Innovation Programme under the project 
\textit{ExaHyPE}, grant agreement number no. 671698 (call FETHPC-1-2014). 

The authors also acknowledge the Leibniz Rechenzentrum (LRZ) in Munich, Germany, for awarding us access to the \textit{SuperMUC} 
supercomputer, as well as the support of the HLRS in Stuttgart, Germany, for awarding access to the \textit{Hazel Hen} supercomputer. 

\reds{The authors would also to thank both the anonymous reviewers for their encouraging comments and remarks that allow us to present a more readable and overall higher-quality version of this paper.}

\bibliographystyle{plain}    
\bibliography{references7}

\end{document}